\documentclass{amsart}
\usepackage{amsmath,amsfonts,amssymb,amscd,multicol}
\usepackage[all]{xypic}
\usepackage[normalem]{ulem}
\usepackage{verbatim,color}
\usepackage{bbm}

\usepackage{fullpage}

\newcommand{\invlim}{\underleftarrow{\lim}}

\DeclareMathOperator{\Ext}{Ext}

\DeclareMathOperator{\Pic}{Pic}

\DeclareMathOperator{\ann}{ann}

\newcommand{\on}{\operatorname}

\DeclareMathOperator{\rProj}{Proj}

\DeclareMathOperator{\Hom}{Hom}

\newcommand{\gldim}{\on{gl.dim}}
\newcommand{\GK}{\on{GKdim}}
\newcommand{\GKdim}{\on{GKdim}}
\newcommand{\rGr}{\on{Gr}}

\newcommand{\rproj}{\on{proj}}

\newcommand{\rgr}{\on{gr}}
\newcommand{\rfd}{\on{tors}}

\newcommand{\mc}{\mathcal}
\newcommand{\mb}{\mathbb}

\newcommand{\mf}{\mathfrak}

\newcommand{\kk}{\Bbbk}
\newcommand{\Proj}{\operatorname{Proj}}

\numberwithin{equation}{section}
 \theoremstyle{plain}
\newtheorem{theorem}[equation]{Theorem}

\newtheorem{proposition}[equation]{Proposition}

\theoremstyle{definition}
\newtheorem{question}[equation]{Question}
\newtheorem{definition}[equation]{Definition}

\newtheorem{standing-hypothesis}[equation]{Standing Hypothesis}
\newtheorem{example}[equation]{Example}

\title{Artin-Schelter Regular Algebras}

\author{Daniel Rogalski}

\dedicatory{To S. Paul Smith on the occasion of his 65th Birthday}

\begin{document}
\maketitle
\tableofcontents

\section{Introduction}
\label{sec:intro}
There are two aspects to this survey.  First, in Sections 1-3 we give an introduction to the main ideas and fundamental results in the theory of Artin-Schelter regular algebras, suitable for readers new to the subject who have some basic background in noncommutative algebra, homological algebra and algebraic geometry.  For those looking for more details, in particular for beginning graduate students, we recommend starting with the course notes in \cite{Rog16} instead.  Second, in the remainder of the paper we give an overview of the most important themes and results of the past 30 years of research in the subject. This literature review, mostly designed to point readers to the original papers where they can find more information, will not be able to linger long on any one topic or idea.  It could not possibly be exhaustive given the vast literature in this area.  We apologize in advance to the authors of papers we do not discuss.   

Throughout these notes, let $\kk$ be a base field, which for convenience we assume is algebraically closed, though some results do not require this assumption.  A polynomial algebra in finitely many variables over $\kk$, say $A =   \kk[x_1, \dots, x_m]$, is one of the most basic examples in commutative ring theory.  Noncommutative ring theory has long been influenced by the beautiful and strong results about commutative rings, motivated by number theory and algebraic geometry, that formed part of the core of the development of modern algebra.  Because polynomial rings are so important, it is natural to ask what noncommutative algebras over a field play a similar role in the noncommutative theory as commutative polynomial algebras do in the commutative theory.

As phrased, this question is vague, since any answer depends on what aspect of the polynomial ring one wants to generalize.  For example, a polynomial ring $\kk[x_1, \dots, x_m]$ is a free object on a set with $m$ elements in the category of commutative $\kk$-algebras.  The noncommutative analog in this sense is the free object on a set with $m$ elements in the category of all associative $\kk$-algebras, namely, the free associative algebra $\kk \langle x_1, \dots, x_m \rangle$.  However, free algebras in more than one variable are unwieldy and have very different ring-theoretic properties from a commutative polynomial ring.  For example, such algebras are not noetherian.  What we seek is a class of rings which are the same size as polynomial rings in some sense, and which have similarly good properties.

One reasonable answer to the question, which goes back to the earliest study of noncommutative rings, is the notion of an iterated Ore extension. Note that we can think of (or define) a polynomial ring via iterating a one variable extension:
\[
\kk[x_1, \dots, x_m] = ( \dots ((\kk[x_1])[x_2]) \dots )[x_m].
\]
Ore defined a noncommutative analog of a single variable extension.
\begin{definition}[{\O}ystein Ore; 1933]
Let $R$ be a ring, $\sigma: R \to R$ an automorphism, and $\delta: R \to R$ a
$\sigma$-derivation, that is, a linear map satisfying $\delta(ab) = \delta(a) b+ \sigma(a) \delta(b)$ for all $a, b \in R$.  The \emph{Ore extension} $R[x; \sigma, \delta]$ is the left free $R$-module $R \oplus Rx \oplus Rx^2 \oplus \dots$ with a new multiplication induced by the rule $xr = \sigma(r) x + \delta(r)$ for all $r \in R$.
\end{definition}
\noindent The reader can find more details on this construction in \cite[Chapter 2]{GW04}.
%It is elementary to check that there is a unique associative ring structure on $R[x; \sigma, \delta]$ determined by the rule $xr = \sigma(r) x + \delta(r)$, which tells you how to move an element in $R$ from the right side of $x$ to the left.  %There is also a converse, due to Cohn, which says that if $S$ is a ring containing $R$ and there is $x \in S$ such that $\{1, %x, x^2 \dots \}$ is a basis for $S$ as a free left and right $R$-module, then $S$ is an Ore extension $R[x; \sigma, \delta]$ %for some $\sigma$ and $\delta$.  
Since an Ore extension is a noncommutative polynomial extension, \emph{iterated Ore extensions} of the form
\begin{equation}
\label{eq:iterated}
\kk[x_1; \sigma_1, \delta_1][x_2; \sigma_2, \delta_2] \dots [x_m; \sigma_m, \delta_m],
\end{equation}
where all of the $\sigma_i$ and $\delta_i$ are $\kk$-linear, give a class of noncommutative $\kk$-algebras which can be thought of as noncommutative analogs of commutative polynomial rings over $\kk$.  These rings give one possible answer to our original question.  Indeed, these rings have many good properties, for example they are noetherian by a noncommutative version of the Hilbert basis theorem.  In addition, \eqref{eq:iterated} has the same size as a polynomial algebra in $m$ variables, in the sense that the set $\{ x_1^{i_1} \dots x_m^{i_m} | i_1, \dots, i_m \geq 0 \}$ is a $\kk$-basis for the algebra.  Many important examples, such as certain quantum groups, can be realized as iterated Ore extensions.  There are many interesting and deep questions about them, such as how to understand the structure of their prime spectra.

In the late 1980's, Artin and Schelter introduced a different class of noncommutative analogs of polynomial rings, motivated instead by the \emph{homological} properties polynomial rings have as graded rings. To study some of these examples, completely new techniques coming from algebraic geometry were required.  This was the seed of an entire subject of noncommutative projective algebraic geometry that has flowered ever since.

To describe the class of algebras Artin and Schelter studied, let us first recall the basic definitions of graded rings. 
\begin{definition}
A $\kk$-algebra is \emph{$\mb{N}$-graded} if there is a $\kk$-space decomposition $A = A_0 \oplus A_1 \oplus \ldots \oplus A_n \oplus \dots $ such that $A_i A_j \subseteq A_{i+j}$ for all $i,j$.
If in addition $A_0 = \kk$, then $A$ is \emph{connected}.  
\end{definition}
In this paper our primary interest is \emph{finitely generated} connected $\mb{N}$-graded $\kk$-algebras.  Note that a commutative polynomial ring $A = \kk[x_1, \dots, x_m]$ has various $\mb{N}$-gradings, by
assigning arbitrary degrees $\deg x_i = d_i$ and letting $A_n$ be the $\kk$-span of 
$\{ x_1^{i_1} x_2^{i_2} \dots x_m^{i_m} | d_1 i_1 + d_2 i_2 + \dots + d_m i_m = n \}$.
In fact, such graded polynomial rings can be characterized as the  commutative graded algebras of finite global dimension. The following result was proved by Smith and Zhang (unpublished), see \cite[Remark after Corollary 3.2]{StZh97}.
\begin{theorem}
\label{thm:comm-reg}
A connected graded finitely generated commutative $\kk$-algebra $A$ with $\on{gl.dim}(A) = m$ is isomorphic to $\kk[x_1, \dots, x_m]$ with some grading.
\end{theorem}

For a connected graded algebra, let $A_{\geq 1} = \bigoplus_{n \geq 1} A_n$ and write $\kk$ for the  \emph{trivial module} $\kk = A/A_{\geq 1}$.  This is both a left and right $A$-module, and we write $_A \kk$ for its left structure and $\kk_A$ for the right structure.  
%In general we work with  left modules unless otherwise noted.

\begin{definition}[Artin-Schelter, 1987]
\label{def:regular}
Let $A$ be a connected finitely generated graded $\kk$-algebra.  The algebra $A$ is \emph{regular} if
\begin{enumerate}
\item $A$ has finite global dimension $d < \infty$;
\item $A$ has polynomial growth, i.e. $f(n) = \dim_{\kk} A_n$ is bounded above by a polynomial function of $n$; and
\item $A$ is \emph{Gorenstein}, i.e.  $\on{Ext}^i_A(_A \kk, A) \cong \begin{cases}  0 & i \neq d \\ \kk_A &  i = d \end{cases}$.
\end{enumerate}
\end{definition}
Algebras satisfying Definition~\ref{def:regular} are now called Artin-Schelter regular or AS regular.  There are other noncommutative notions of regular ring, but in this paper if we just say ``regular" we will mean ``AS regular". 
Also, while it is not required in the definition that $A$ be \emph{generated in degree $1$}, that is, generated as a $\kk$-algebra by the elements of $A_1$, all of the initial results in the subject have focused on that case, and we will tacitly assume it unless we say otherwise.

There is a related more general notion:  $A$ is called \emph{AS Gorenstein} if it satisfies (2) and (3) of the definition above but 
(1) is replaced by the weaker condition that $A$ has finite injective dimension as an $A$-module.  AS Gorenstein rings will only be mentioned in passing in this paper.

While commutative graded rings of finite global dimension are polynomial rings, by Theorem~\ref{thm:comm-reg}, noncommutative graded rings of finite global dimension do not behave like analogs of polynomial rings without further restrictions.  For example, a free associative algebra $\kk \langle x_1, \dots, x_m \rangle$ has global dimension 1, regardless of $m$, and is not noetherian when $m \geq 2$.  But this example satisfies neither (2) nor (3) in the definition of AS regular when $m \geq 2$.

Here are some other examples that show the importance of conditions (2) and (3) in the definition of regularity.  We don't justify the claims made here, but they all follow using methods similar to those reviewed in the next section (but easier).

\begin{example}
\label{ex:bad}
Let $A = \kk\langle x, y \rangle/(yx)$.  Then $A$ has a $\kk$-basis consisting of words $\{ x^i y^j | i, j \geq 0\}$, the same as for a polynomial ring in two variables.  In particular, $A$ has polynomial growth. Also, one may show that $\on{gl.dim}(A) = 2$.  However $A$ is patently not a domain, and it is easy to show it is not noetherian.  One may also see that $A$ is not Gorenstein, that is, condition (3) fails.
\end{example}

\begin{example}
Let $A = \kk \langle x, y, z \rangle/(x^2 + y^2 + z^2)$.  Then $A$ is a graded domain with
global dimension $2$ and the Gorenstein property (3), but with $f(n) = \dim_{\kk} A_n$
having exponential growth.  Again, $A$ is non-noetherian.
\end{example}

On the other hand, regular algebras of global dimension at most $2$ are easy to classify, and do have good properties similar to those of commutative polynomial rings. 

\begin{example}
If $A$ is a $\kk$-algebra which is regular of global dimension $1$, then $A \cong \kk[x]$.
\end{example}

\begin{example}
\label{ex:dim2}
Let $A$ be a regular algebra of global dimension $2$.  Then either $A \cong A_q = \kk \langle x, y \rangle/(yx - qxy)$ for some $0 \neq q \in \kk$ (the \emph{quantum plane}), or $A \cong A_J = \kk \langle x, y \rangle/(yx -xy - x^2)$ (the \emph{Jordan plane}).  
\end{example}

If $R$ is an $\mb{N}$-graded ring, a general Ore extension $S = R[x; \sigma, \delta]$ need not have a natural grading.  
However, if $\sigma$ is a graded automorphism ($\sigma(R_n) \subseteq R_n$ for all $n$) and $\delta$ is a $\sigma$-derivation of degree $1$ ($\delta(R_n) \subseteq R_{n+1}$ for all $n$), then $S$ is again graded, with $S_n = \sum_{i+j = n} R_i x^j$.  One may check that in fact $A_q \cong \kk[x][y; \sigma, \delta]$ as a graded ring, where $\sigma(x) = qx$ and $\delta(x) = 0$, and similarly $A_J \cong \kk[x][y; \sigma, \delta]$, where $\sigma(x) = x$ and $\delta(x) = x^2$.  It is known that a graded iterated Ore extension is indeed regular, but the notion of regular algebra turns out to be far more general.  Already in dimension $3$ there are interesting examples of regular algebras which cannot be constructed as iterated Ore extensions.

\textbf{Acknowledgments:}  We thank Paul Smith and Michaela Vancliff for  helpful comments and suggestions.  We especially thank the referee for their thoughtful advice which led to many significant improvements to this survey.  We also thank our many mentors and friends with more expertise than ourselves, including Paul Smith to whom this article is dedicated, for conversations over the years from which we learned about many of the beautiful results of this subject.  Finally, we thank Michaela Vancliff for writing the companion article \cite{Van23} in response to a comment from the referee.

\section{Regular algebras of dimension 3}
\label{sec:reg3}

The classification of regular algebras of dimension 3 was the goal of the earliest work in the subject, starting with the original paper by Artin and Schelter, and continuing in two papers of Artin, Tate, and Van den Bergh.  In this section we review the main ideas of these papers, which have influenced all subsequent work.   

\subsection{The work of Artin and Schelter}
\label{sec:AS}

We would like to give the reader a sense of why the axioms of an AS regular algebra are restrictive enough to give any hope of classifying those of low dimension.  In fact, we will see that all of the properties in the definition of regular can be encoded in terms of the properties of the free resolution of the trivial module $\kk$.  The reference for this entire subsection is the paper \cite{AS87}.

Let us first make a comment on morphisms between free modules.  An $A$-module morphism between free left modules of finite rank, say  $f: A^{\oplus r} \to A^{\oplus s}$, can be represented by an $r \times s$ matrix $M$ with entries in $A$, where we represent the elements of the free modules by row vectors and $f$ by right multiplication by the matrix $M$. 
Similarly, a morphism between free right modules $g: A^{\oplus s} \to A^{\oplus r}$ can also be represented by an $r \times s$-matrix $M$, where we now represent the free modules as columns and $g$ by left multiplication by the matrix.  The contravariant functor $\Hom_A(-, A)$ interchanges the categories of left and right free modules of finite rank, and leaves the matrices representing maps the same if we use these conventions.

Let $A$ be a connected graded finitely generated $\kk$-algebra.  We work in the category of $\mb{Z}$-graded left $A$-modules $M = \bigoplus_{n \in \mb{Z}} M_n$.  For any such module $M$ and $i \in \mb{Z}$, let $M(i)$ be the graded module with the same ungraded module structure $M$, but with new grading $M(i)_j = M_{i+j}$.  Every finitely generated graded module $M$ has a graded free resolution 
\begin{equation}
\label{eq:res}
 \to F_n \to \dots \overset{\partial_1}{\to} F_1 \overset{\partial_0}{\to} F_0 \to M \to 0,
\end{equation}
where each $F_i$ is a direct sum of various graded free modules $A(j)$ and the maps preserve degree.  Let us assume 
that each $F_i$ has finite rank, which we can always achieve if $A$ is noetherian.  A resolution such as \eqref{eq:res} is \emph{minimal} if the matrices representing the maps $\partial_i$ have all of their entries in $A_{\geq 1}$.  The minimal resolution is unique up to isomorphism of complexes, and as the word minimal suggests, it is the resolution for which the ranks of the free modules $F_i$ are smallest.  

For any connected graded finitely generated $\kk$-algebra $A$, it is well known that the global dimension of the algebra $A$ is the same as the projective dimension of the trivial module $\kk$ (as a left or right module).  Thus $\gldim(A) = d$ if and only if $\on{proj.dim}(_A \kk) = d$, which is if and only if the minimal resolution of $\kk$ as in \eqref{eq:res} has $F_i = 0$ for all $i \geq d+1$.

Now let $A$ be regular of dimension $3$, and generated in degree $1$.  Let $A$ have a minimal presentation $A \cong \kk \langle x_1, \dots, x_n \rangle/(r_1, \dots, r_m)$, where $\deg(x_i) = 1$ and each $r_i$ is a homogeneous relation with $\deg(r_i) = d_i$, say.  The degrees of the generators and relations determine the ranks and shifts of the first few graded free modules in the minimal resolution of ${}_A \kk$.  Specifically, since $\on{proj.dim}(k) = 3$, the minimal graded free resolution of the trivial left module $\kk = A/A_{\geq 1}$ has the form
\[
0 \to  F \to \bigoplus_{i=1}^m A(-d_i) \to \bigoplus_{i=1}^n A(-1) \to A \to \kk \to 0,
\]
for some graded free module $F$.

Next, consider the Gorenstein condition.  In the definition of regular we ignored shifts, but 
$\Ext^i_A(\kk, A)$ is naturally a $\mb{Z}$-graded right $A$-module, and so we can actually assume that 
$\Ext^3_A(\kk, A) \cong \kk_A(\ell)$ as graded right $A$-modules, for some shift of grading $\ell$.
Calculating $\Ext$ by applying $\Hom_A(-, A)$ to the (deleted) free resolution above, $\Ext^i(\kk, A)$ is the $i$th 
cohomology group of a complex
\[
0 \to A \to \bigoplus_{i=1}^n A(1) \to \bigoplus_{i=1}^m A(d_i) \to \Hom(F, A) \to 0.
\]
Since $\Ext^i(\kk, A) = 0$ for $i \neq 3$, the complex is exact in every spot except at $\Hom(F, A)$.  We also know $\Ext^3(\kk, A) = \kk(\ell)$, so there is an exact sequence of right modules 
\[
0  \to A \to \bigoplus_{i=1}^n A(1) \to \bigoplus_{i=1}^m A(d_i) \to \Hom(F, A)  \to \kk(\ell) \to 0.
\]
Shifting all terms by $-\ell$ we get 
\[
0  \to A(-\ell) \to \bigoplus_{i=1}^n A(1-\ell) \to \bigoplus_{i=1}^m A(d_i-\ell) \to \Hom(F, A)(-\ell)  \to \kk \to 0.
\]
This must be a minimal graded free resolution of the right trivial module $\kk_A$, since the matrices representing the maps have 
not changed, so they still have entries in $A_{\geq 1}$.  Also, the ranks and shifts of its first few terms are again determined by the number of generators and the degrees of the relations, so we conclude that $F \cong A(-\ell)$; $d_i- \ell = -1$ for all $i$; and the minimal number $m$ of relations is the same as the minimal number $n$ of generators of $A$.  The original minimal resolution of ${}_A \kk$ now must look like
\[
0 \to  A(-\ell) \to \bigoplus_{i=1}^m A(-\ell+1) \to \bigoplus_{i=1}^m A(-1) \to A \to \kk \to 0.
\]
To summarize, the Gorenstein condition implies that the minimal resolutions of the left and right trivial modules are interchanged (up to shift of grading) by the operation $\Hom_A(-, A)$.  This forces the minimal resolution to have a rather symmetric form.

Finally, we can relate the form of the minimal free resolution of $\kk$ to the growth of the graded pieces of the algebra.
The \emph{Hilbert series} of a connected finitely generated $\mb{N}$-graded algebra is the power series $h_A(t) = \sum_{n = 0}^{\infty} (\dim_{\kk} A_n) t^n$, which can be manipulated in the power series ring $\mb{Z}[[t]]$.  Using the form of the minimal free resolution of $\kk$, it is straightforward to see that 
\[
h_A(t) = \sum_{n =0}^{\infty} (\dim_{\kk} A_n) t^n = \frac{1}{(1 - mt + mt^{\ell} - t^{\ell+1})}.
\]
Similarly, the Hilbert series of any connected finitely graded algebra of finite global dimension is of the form $\frac{1}{p(t)}$, 
where the coefficients of the polynomial $p(t)$ can be determined in a simple way from the ranks and shifts of the graded free modules in the minimal free resolution of $\kk$.  Some basic complex analysis can be used to show that if $p(t)$ is a monic polynomial in $\mb{Z}[t]$ and $\frac{1}{p(t)} = 1 + a_1 t + a_2 t^2 + \dots$ in $\mb{Z}[[t]]$, then $f(n) = a_n$ is a polynomially bounded function of $n$ if and only if all zeroes of $p(t)$ in $\mb{C}$ are roots of unity \cite[Proposition 2.14]{ATV2}.  Thus, the polynomial growth condition (2) in the definition of regularity forces $1 - mt + mt^{\ell} - t^{\ell+1}$ to have only roots of unity for zeroes.  A little calculus shows this happens if and only if either $m =3, \ell = 2$ or $m = 2, \ell = 3$.  In conclusion, either $A$ has $3$ generators and $3$ quadratic relations, with $h_A(t) = \frac{1}{(1-t)^3}$, or $A$ has $2$ generators and $2$ cubic relations, with $h_A(t) = \frac{1}{(1-t)^2(1-t^2)}$.  These are the same as the Hilbert series for a polynomial ring $\kk[x, y,z]$ with variables of degrees $1, 1, 1$ or $1, 1, 2$, respectively.  (Note, however, that the polynomial ring $\kk[x, y, z]$ with generators of degrees $1, 1, 2$ is not generated in degree $1$, and so is not an example of the type considered in this section.)  
Based on the degrees of their relations, the two kinds of regular algebras of dimension $3$ are referred to as quadratic or cubic.

Both kinds of Hilbert series were treated uniformly by Artin and Schelter, but here let us focus on the quadratic case, where $A = \kk \langle x_1, x_2, x_3 \rangle/(r_1, r_2, r_3)$ has $3$ minimal generators and $3$ minimal relations of degree $2$.  Let us consider the maps in the free resolution
\[
0 \to  A(-\ell) \overset{\partial_3}{\to} \bigoplus_{i=1}^m A(-\ell+1) \overset{\partial_2}{\to} \bigoplus_{i=1}^m A(-1) \overset{\partial_1}{\to} A \to \kk \to 0.
\]
The map $\partial_i$ is given by right multiplication by some matrix $M_i$, and we identify $\partial_i$ and $M_i$.
It is standard that we can take $\partial_1$ to be the $3 \times 1$ matrix $(x_1, x_2, x_3)^t$.  Because of the symmetry coming from the Gorenstein condition, the three entries of $\partial_3$ also must be a basis of $A_1$.  By making a change of basis 
in the $\bigoplus_{i=1}^m A(-\ell+1)$-spot, we can assume that $\partial_3 = (x_1, x_2, x_3)$.  This also changes the matrix representing $\partial_2$, but no assumptions about it have been made.  The free resolution now looks like  
\[
0 \to  A(-3) \overset{(x_1, x_2, x_3)}{\to} \bigoplus_{i=1}^3 A(-2) \overset{M}{\to} \bigoplus_{i=1}^3 A(-1) \overset{(x_1, x_2, x_3)^t}{\to} A \to \kk \to 0.
\]
Note that the entries of the $3 \times 3$ matrix $M$ are in $\kk x_1 + \kk x_2 + \kk x_3$.  By the standard way the resolution of $\kk$ begins, when calculated in the free algebra $\kk \langle x_1, x_2, x_3 \rangle$ the product $M (x_1, x_2, x_3)^t$ is equal to $(r_1, r_2, r_3)^t$ for the three minimal relations $r_1, r_2, r_3$.  Again by the symmetry coming from the Gorenstein condition, the product of the final two matrices (taken in $\kk \langle x_1, x_2, x_3 \rangle$) is $(x_1, x_2, x_3) M = (s_1, s_2, s_3)$, where $s_1, s_2, s_3$ is also a basis for the span of the relations $\kk r_1 + \kk r_2 + \kk r_3$.  So there is a matrix $Q \in \on{GL}_3(\kk)$ with $(x_1, x_2, x_3) M = M (x_1, x_2, x_3)^t Q$.

Consider the product $w= (x_1, x_2, x_3) M (x_1, x_2, x_3)^t$, again taken in the free algebra $\kk \langle x_1, x_2, x_3 \rangle$.  Note that $w$ determines the relations $r_1, r_2, r_3$, since $w$ has a unique expression $w = x_1 r_1 + x_2 r_2 + x_3 r_3$ in the free algebra.  Using the equation $(x_1, x_2, x_3) M = M (x_1, x_2, x_3)^t Q$, one checks that there is a graded automorphism $\sigma$ of $\kk \langle x_1, x_2, x_3 \rangle$ (given by the matrix $Q$) such that under the linear map $\tau$ which operates on degree $3$ monomials by  $x_{i_1}x_{i_2}x_{i_3} \to \sigma(x_{i_3})x_{i_1}x_{i_2}$, one has $\tau(w) = w$.  
%In fact $\sigma$ preserves the ideal of relations and so descends to an automorphism $\sigma$ of $A$ called the \emph{Nakayama automorphism}.  
This says that $w$ is a \emph{twisted superpotential}, in more modern terminology (not used by Artin and Schelter) we explain later in Section~\ref{sec:CYtoAS}.

Using some simple representation theory, Artin and Schelter wrote down a list of possible twisted superpotentials corresponding to diagonalizable $Q$'s.  They were able to show certain superpotentials led to regular algebras by direct calculation, namely 
the theory of noncommutative Gr\"obner bases, which played a large role in this early work.  Schelter's Affine computer 
algebra package was a key tool.  For the reader unfamiliar with noncommutative Gr\"obner bases, we refer to Bergman's survey \cite{Berg78} for a general introduction, or see the author's course notes in \cite[Section 1]{Rog16}.  In particular, the Gr\"obner basis method allows one to compute the Hilbert series of a finitely presented algebra in some cases, and can also aid in other calculations, such as proving that elements are non-zero-divisors.  Such information is helpful in proving that a proposed free resolution of $\kk$ is actually exact.

Some algebro-geometric analysis showed that, up to isomorphism, all regular algebras of dimension $3$ whose twisted superpotentials have a diagonalizable $Q$ belong to one of a finite number of ``types", where each type is a parametrized family of superpotentials.  Conversely, each type gives a regular algebra for generic choices of the parameters.  The types of quadratic algebras were named $A, B, H, S_1, S_1', S_2$.  Let us give an example of one of the simpler types.
\begin{example}
\label{ex:S2}
Type $S_2$ is given by the twisted superpotential 
\[
w = x_1 x_3 x_1 + \alpha x_3 x_1^2 + \alpha^{-1} x_1^2 x_3 + x_2 x_3 x_2 - \alpha x_3 x_2^2 - \alpha^{-1} x_2^2 x_3
\]
for any $0 \neq \alpha \in \kk$.  This is a twisted superpotential where $\sigma(x_1) = \alpha^{-1} x_1$, $\sigma(x_2) = -\alpha^{-1} x_2$, $\sigma(x_3) = \alpha^2 x_3$.  Since $w = x_1 r_1 + x_2 r_2 + x_3 r_3$ for the corresponding relations $r_1, r_2, r_3$, 
we have a family of algebras $A(\alpha)$ with parametrized relations:
\[
A(\alpha) = \kk \langle x_1, x_2, x_3 \rangle/(\alpha x_3x_1 +  x_1 x_3, \ \alpha x_3x_2- x_2x_3, \ x_1^2 - x_2^2).
\]
\end{example}
The family $S_2$ is quite simple, because it depends on only one parameter and each relation involves only two of the variables.
In fact, an algebra in this family is easily represented as an iterated Ore extension $\kk[x_1+x_2][x_1-x_2; \sigma][x_3;\tau]$, and one can use this to see that $A(\alpha)$ is regular for all nonzero $\alpha$.  On the other hand, for some of the other families in the classification, the question of which values of the parameters give a regular algebra turned out to be quite complicated, and was not answered until the later work of Artin, Tate, and Van den Bergh in \cite{ATV1}.

\subsection{Linear modules}

The following type of module plays a crucial role in the theory of AS regular algebras.
\begin{definition}
Let $A$ be a connected $\mb{N}$-graded $\kk$-algebra, generated in degree $1$.
A cyclic graded $A$-module $M$, generated in degree $0$, is called \emph{linear} if it has Hilbert series $h_M(t) = \bigoplus_{n = 0}^{\infty} \dim_{\kk}(M_n) t^n = \frac{1}{(1-t)^d}$ for some $d \geq 0$.
If $d = 1$, $M$ is called a \emph{point module}, and if $d = 2$, $M$ is called a \emph{line module}.
\end{definition}
Point modules were essential to Artin, Tate, and Van den Bergh's approach, which we review in the next section, to the classification of AS regular algebras of dimension $3$. To help the reader picture such modules, we note that since $\frac{1}{(1-t)} = 1 + t + t^2 + t^3 + \dots$, a point module $P$ has $\dim_{\kk} P_n = 1$ for all $n \geq 0$.    The polynomial ring $\kk[x]$ is a point module over itself, and over any connected graded algebra that maps onto $\kk[x]$.  Similarly, a line module $L$ has $\dim_{\kk} L_n = n+1$ for all $n \geq 0$, and $\kk[x, y]$ is a line module over any connected graded algebra mapping onto it.

If $A = \kk[x_1, \dots, x_n]$ is a commutative polynomial ring, then one may prove that a linear module $M$ with $h_M(t) = \frac{1}{(1-t)^d}$ is isomorphic to $A/I$ for $I = \ann_A(M)$, where $I$ is generated as an ideal by any $\kk$-subspace $V$ of $A_1$ with $\dim_{\kk} V = n-d$.  Thus such linear modules over $A$ are in bijective correspondence with the projective linear subspaces of $\Proj A \cong \mb{P}^{n-1}$ that are isomorphic to $\mb{P}^{d-1}$.  In particular, the isomorphism classes of point modules over $A$ are simply in bijection with the $\kk$-points of $\mb{P}^{n-1}$. More generally, when $A = \kk[x_1, \dots, x_n]$ the isomorphism classes of point modules over any factor ring $A/I$ by a homogeneous ideal $I$ are in bijection with the $\kk$-points of $\Proj A/I$.

Now let $A = \kk \langle x_1,\dots, x_n \rangle/(r_1, \dots, r_m)$ be a finitely presented connected graded $\kk$-algebra.  For any $d$, define a \emph{truncated point module of length $d+1$} to be a graded $A$-module, generated in degree $0$, with Hilbert series $1 + t + \dots + t^d$ (length here refers to the composition length).  Standard arguments show that the isomorphism classes of truncated point modules of length $d+1$ for $A$ are parametrized by a closed subscheme $X_d$ of $(\mb{P}^{n-1})^{\times d}$.  In fact, the scheme $X_d$ can be easily calculated as the common zero locus in $(\mb{P}^{n-1})^{\times d}$ of multilinearizations of the relations $r_i$.  We will speak of parametrization informally in these notes.  It means that these truncated point modules, up to isomorphism, are in bijection with the $\kk$-points of the scheme $X_d$, in a natural way which is compatible with base ring extension; see \cite[Section 3]{ATV1}.
%To make the notion precise, one shows that when one replaces the field $\kk$ with an arbitrary commutative noetherian base ring $R$, the $R$-points of $\Gamma_d$ correspond to the isomorphism classes of graded modules $M$ over $A \otimes_{\kk} R$, generated in degree $0$, with $M_0 = R$ and $M_i$ locally free of rank $1$ over $R$ for each $1 \leq i \leq d$.  $\Gamma_d$ is what is known as a fine moduli space in algebraic geometry. 

For each $d \geq 1$ there is a morphism of schemes $\phi_d: X_d \to X_{d-1}$, which corresponds to taking a truncated point module $M$ of length $d+1$ and truncating it further to get a truncated point module $M_{\leq (d-1)}$ of length $d$.  Then the inverse limit of schemes $X = \invlim \, X_d$ parametrizes the isomorphism classes of full point modules for $A$.  In specific examples, one often sees through direct calculation that this inverse system stabilizes, i.e. that there is $\ell \geq 1$ such that $\phi_d$ is an isomorphism for all $d > \ell$.  In this case $X = X_{\ell}$ is a projective scheme which parametrizes the point modules for $A$, and is called the \emph{point scheme} for $A$.  If the inverse limit should happen not to stabilize, then we can still treat $X$ as an inverse limit of schemes, but it may not be an object in the category of schemes itself.  If $M$ is  point module for $A$, then we can lop off the zero-degree piece of $M$ and shift everything over, obtaining another point module $M_{\geq 1}(1)$.    This operation translates into a map $\sigma: X \to X$ of the parametrizing space $X$.  The map $\sigma$ is a morphism of schemes in the case that $X$ is a scheme, and is typically an automorphism of $X$.

Here are a few simple examples of the theory above.
\begin{example}
If $A = A_q = \kk \langle x, y \rangle/(yx- qxy)$ for $q \neq 0$, then
the (left) point modules up to isomorphism are the modules $A/A (\alpha x + \beta y)$, parameterized by points $(\alpha: \beta) \in \mb{P}^1$.  The induced map $\sigma: \mb{P}^1 \to \mb{P}^1$ is given by $\sigma(\alpha: \beta) = (q \alpha: \beta)$.
\end{example}

\begin{example}
\label{ex:S2revisited}
As in Example~\ref{ex:S2}, consider the family $S_2$ of regular algebras of dimension $3$, given by the presentation 
\[
A(\alpha) = \kk \langle x_1, x_2, x_3 \rangle/(\alpha x_3x_1 +  x_1 x_3, \ \alpha x_3x_2- x_2x_3, \ x_1^2 - x_2^2)
\]
for nonzero $\alpha$.  The scheme $X_2$, which parametrizes truncated point modules of length $3$, is the subscheme of $\mb{P}^2 \times \mb{P}^2$ defined by the vanishing of the multilinearized relations, that is 
\begin{gather*}
\{ ((x_1: x_2: x_3), (y_1: y_2: y_3)) \, | \, \alpha x_3y_1 + x_1y_3 = \alpha x_3y_2-x_2y_3 = x_1y_1-x_2y_2 = 0 \} \\ 
= X_2 \subseteq \mb{P}^2 \times \mb{P}^2. 
\end{gather*}
One can check that the projection of $X_2$ onto either coordinate is the subscheme $E$ of $\mb{P}^2$ given by the union of 
the three coordinate lines, and $X_2 = \{ (\sigma(p), p) | p \in E \}$ for an automorphism $\sigma: E \to E$ which interchanges two of the lines.  It is known that the inverse system already stabilizes at this point, so that $X = X_2 \cong E$ is the point scheme, and the truncation shift on point modules corresponds to the automorphism $\sigma$.
\end{example}

For regular algebras of dimension $3$, the only interesting linear modules are the point modules.  For example, if $A$ is regular of dimension $3$ with Hilbert series $\frac{1}{(1-t)^3}$, then $A$ is known to be a domain, and it easily follows that the line modules of $A$ are simply given by the modules $A/Ay$ for any nonzero $y \in A_1$.  In particular, the line modules are parametrized by $\mb{P}^2$.  Line modules, and other higher-dimensional linear modules, become more relevant in the investigation of regular algebras of larger dimension.

\subsection{The work of Artin-Tate-Van den Bergh}

Artin, Tate and Van den Bergh gave a geometric description of regular algebras of dimension $3$, which refined 
the original classification results of Artin and Schelter.  We state their results only in the quadratic case, for simplicity.  However, their work also fully handles the cubic case, in a completely analogous manner.

\begin{theorem} \cite[Theorems 1,2]{ATV1}
\label{thm:ATV1}
Let $A$ be an AS regular algebra of dimension $3$ with three degree one generators and three quadratic relations. 
\begin{enumerate}
    \item The point modules for $A$ are parametrized by either $E = \mb{P}^2$ or a degree three divisor $E \subseteq \mb{P}^2$, and translation-shift gives an automorphism $\sigma: E \to E$; 
    \item Defining the sheaf $\mc{L} = \mc{O}_{\mb{P}^2}(1) \vert_E$ with corresponding element $\lambda \in \Pic E$, then $(\sigma -1)^2 \lambda = 0$ always, while $(\sigma -1) \lambda = 0$ if and only if $E = \mb{P}^2$; 
    \item $B = B(E, \mc{L}, \sigma) = \bigoplus_{n \geq 0} H^0(E, \mc{L} \otimes \sigma^*(\mc{L}) \otimes \dots \otimes (\sigma^{n-1})^*(\mc{L}))$ has a natural ring structure, and there is a normal element $g \in A_3$ such that $A/gA \cong B$. When $E = \mb{P}^2$ then $g = 0$; otherwise, $g$ is a non-zero-divisor in $A$.
\end{enumerate}
\end{theorem}
The ring $B$ appearing in part (3) is called a \emph{twisted homogeneous coordinate ring}.  The multiplication in the ring is defined for $x \in B_m$, $y \in B_n$, by $x \otimes (\sigma^m)^*(y)$; here, for $y \in H^0(\mc{L} \otimes \sigma^*(\mc{L}) \otimes \dots \otimes (\sigma^{n-1})^*(\mc{L}))$ we note that pulling back by $\sigma^m$ yields an element in $H^0((\sigma^m)^*\mc{L} \otimes \dots \otimes (\sigma^{m+n-1})^*(\mc{L}))$, so that $x \otimes (\sigma^m)^*(y)$ gives an element in $H^0(\mc{L} \otimes \sigma^*(\mc{L}) \otimes \dots \otimes (\sigma^{n+m-1})^*(\mc{L}))$, that is, an element of $B_{m+n}$.
Twisted homogeneous coordinate rings $B(X, \mc{L}, \sigma)$ can be defined for any projective scheme $X$, invertible sheaf $\mc{L}$, and automorphism $\sigma$, and are interesting rings in their own right, apart from their applications to regular algebras \cite{AV90}.  Their ring-theoretic properties can often be understood in terms of the associated geometry.  
There is a robust understanding of when they are noetherian, for example \cite{Kee00}.  The particular 
rings $B(E, \mc{L}, \sigma)$ occurring in part (3) of the theorem are always noetherian.  A general lemma for connected graded rings \cite[Section 8]{ATV1} implies that because $A/gA \cong B$, the ring $A$ is also noetherian.

The converse of Theorem~\ref{thm:ATV1} also holds:
\begin{theorem} 
\label{thm:ATV-converse}
\cite[Theorem 3]{ATV1}
Consider the twisted homogeneous coordinate rings $B(X, \mc{L}, \sigma)$ as  above.
\begin{enumerate}
\item
The algebra $B(\mb{P}^2, \mc{O}_{\mb{P}^2}(1), \sigma)$ is regular of dimension $3$ for any automorphism $\sigma$ of $\mb{P}^2$.
\item For any degree $3$ divisor $E \subseteq \mb{P}^2$ and
automorphism $\sigma: E \to E$ such that the element $\lambda \in \Pic E$ corresponding to
$\mc{L} = \mc{O}_{\mb{P}^2}(1) \vert_E$ satisfies the properties $(\sigma -1)^2 \lambda = 0$ and $(\sigma -1) \lambda \neq 0$, the ring $B = B(E, \mc{L}, \sigma)$ has a minimal presentation with $3$ generators of degree $1$ and minimal relations of degrees $2, 2, 2, 3$.  Removing the degree $3$ relation produces an regular algebra $A$ of dimension $3$ with a normal element $g \in A_3$ such that $A/gA \cong B$.
\end{enumerate}
\end{theorem}
The theorem shows that a regular algebra of dimension $3$ with Hilbert series $\frac{1}{(1-t)^3}$ can be recovered from the geometric data $(E, \mc{L}, \sigma)$ coming from the point scheme, and thus the choices of such triples satisfying Theorem~\ref{thm:ATV1}(1)(2), up to a natural notion of equivalence of triples, are in bijective correspondence with the isomorphism classes of quadratic regular algebras of dimension $3$.

\begin{example}
\label{ex:sklyanin3}
Consider the following family of algebras for $a, b, c \in \kk$:
\[
S(a, b, c) = \kk \langle x_1, x_2, x_3 \rangle/(ax_i^2 + bx_{i+1}x_{i+2} + cx_{i+2}x_{i+1} \, | \, i \in \mb{Z}_3 ),
\]
where the indices are taken modulo $3$.
This family was known as ``type $A$" in Artin and Schelter's original classification. Any member of the family is now known as a  \emph{Sklyanin algebra of dimension $3$}, for reasons we describe in Section~\ref{sec:elliptic} below.  The computational methods used in  \cite{AS87} involving Gr\"obner bases did not succeed to understand this family.  Using the geometric approach of \cite{ATV1}, on the other hand, it is relatively simple to calculate the point modules for $S$, and to use the theorems above to show that an algebra in this family is regular precisely when at least two of the parameters $a, b, c$ are nonzero and the equation $a^3 = b^3 = c^3$ does not hold.  When $S$ is regular, for generic values of $a, b, c$ the point modules $E$ are parametrized by a smooth elliptic curve $E$, and the associated automorphism $\sigma: E \to E$ is translation by a point in the group structure of $E$.
\end{example}
We remark that Iyudu and Shkarin \cite{IySh17}, \cite{IySh21} showed many years later the surprising result that a (necessarily infinite) Gr\"obner basis for the ideal of relations in the Sklyanin algebra can in fact be analyzed effectively, and used to find the Hilbert series of this algebra computationally.  In particular, they recovered in a different way the classification of exactly which parameters $a,b,c$ give a regular algebra. 

Point modules are not the only way to approach the geometry associated to regular algebras of dimension $3$.  Bondal and Polishchuk \cite{BoPo93} gave a different point of view on the classification of quadratic regular algebras of dimension $3$, by considering the notion of the $\mb{Z}$-algebra associated to a helix.  Ultimately, this still involves considering line bundles on cubic divisors in $\mb{P}^2$.  This other approach turns out to be especially interesting when applied in the cubic case, where work of Van den Bergh showed that ``cubic regular $\mb{Z}$-algebras" are a bit more general than cubic regular algebras of dimension $3$ \cite{VdB11}.

Given a parametrized family whose generic member is AS regular, it is also interesting to ask what happens for the special values of the parameters where the algebra is not regular.  In the case of the Sklyanin family in Example~\ref{ex:sklyanin3}, the non-regular members are known as \emph{degenerate} Sklyanin algebras.  Their properties were first studied by Walton \cite{Wal09}, who showed that these algebras $S$ all have the same Hilbert series, where $\dim_{\kk} S_n$ grows exponentially, and that they are neither noetherian nor domains.  She also calculated the truncated point schemes $X_d$ for these algebras, and showed that the inverse limit $\invlim \, X_d$ does not stabilize.  However, there is still a surjection $S \to B$, for a kind of generalized twisted homogeneous coordinate ring $B$ formed from the point scheme data.  Smith \cite{Smi12a} noted that all of these degenerate Sklyanin algebras actually have monomial relations after a change of basis of the generators, and that they are all Zhang twists of each other, in the sense we define below in Section~\ref{sec:twists}.  This helps to explain why the degenerate Sklyanin algebras all have essentially the same properties.  See also \cite{DeL17a}.

\section{Basic properties of regular algebras}

As soon as the classification of regular algebras of dimension $3$ was completed, the obvious question was whether similar ideas could be used to classify regular algebras of dimension $4$. The dimension $4$ case has proven to be much more complicated, and so far has  resisted any sort of general classification.  However, the search for a better understanding of dimension $4$ (and higher-dimensional) regular algebras has led to the discovery of many interesting examples and techniques.  In this section we review some additional background and describe some general questions about regular algebras.  Then in the succeeding sections we describe some of the known results, organized by general theme.

One consequence of the classification of regular algebras of dimension $3$ is that such algebras have good ring-theoretic properties.  For example, we have already remarked that they are noetherian, and they also must be domains \cite[Section 3]{ATV2}.  There are a number of other properties that also hold for all known examples of regular algebras, and thus it is natural to conjecture that they hold in general.  It is true that the vast majority of specific examples that have been examined in detail have global dimension at most 5, so it remains possible that the general higher-dimensional picture is different somehow. 

\subsection{GK-dimension}
\label{sec:GK}

 First let us recall the GK-dimension of a module.  The reader can find more details in \cite{KL00}.
\begin{definition}
Let $A$ be a finitely generated $\kk$-algebra, generated by a finite-dimensional $\kk$-subspace $V$.  Let $M$ be a finitely generated right $A$-module, generated by a finite-dimensional subspace $W$.  The \emph{Gelfand-Kirillov dimension} (GK-dimension) of $M$ is $\GK(M) = \displaystyle \limsup_{n \to \infty} \log_n \dim_{\kk}(W V^n)$ (which one may prove is independent of $V$ and $W$).
\end{definition}
Essentially, $\GK(M) = \alpha$ means that the sequence $f(n) = \dim_{\kk} W V^n$ grows roughly at the same rate as $n^{\alpha}$. When $A$ is a connected, finitely generated $\mb{N}$-graded algebra, and $M$ is a finitely generated $\mb{Z}$-graded $A$-module, then 
\[
\GK(M) = \limsup_{n \to \infty} \log_n d(n), \ \text{where}\ d(n) = \displaystyle \dim_{\kk} \bigoplus_{i = -n}^n M_i = \sum_{i=-n}^n \dim_{\kk} M_i.
\]
In particular, applying this to the ring $A$ itself, it is easy to see that $\GK(A) < \infty$ if and only if the sequence $\dim_{\kk} A_n$ has polynomially bounded growth.  So condition (2) in Definition~\ref{def:regular} is equivalent to $\GK(A) < \infty$.

Now suppose that $A$ is a connected $\mb{N}$-graded algebra with finite global dimension.  As mentioned in the discussion of the global dimension $3$ case, by examining the graded free resolution of the trivial module $\kk$, one sees that the Hilbert series of $A$ has the form $h_A(t) = \frac{1}{p(t)}$ for some monic polynomial $p(t) \in \mb{Z}[t]$.  Moreover, $\GK(A) < \infty$ if and only if $p(t)$ has only roots of unity for zeroes.  If instead $\GK(A) = \infty$, then one can prove in fact that $A$ has \emph{exponential growth}, i.e. that $\dim_{\kk} A_n$ grows asymptotically like some function $\beta^n$ with $\beta > 1$.  Stephenson and Zhang showed in \cite[Theorem 1.2]{StZh97} that a connected graded algebra with exponential growth cannot be noetherian  (this paper is also a good place to find proofs of the other facts mentioned above). Thus, finite GK-dimension is a necessary condition for a connected graded algebra of finite global dimension to be noetherian.  This certainly justifies why the polynomial growth condition is needed in the definition of AS regular, if one is interested in noetherian rings.  
On the other hand, finite GK-dimension is not a \emph{sufficient} condition for an algebra of finite global dimension to be noetherian, as Example~\ref{ex:bad} shows.  All known examples of regular algebras are noetherian, so the Gorenstein condition seems to be absolutely crucial in order to guarantee good ring-theoretic properties, in a way that is still not well-understood.  

We mention that in contrast to the results above for algebras of finite global dimension, it is unknown if there exists a connected graded noetherian algebra $A$ (necessarily of infinite global dimension) with subexponential growth, yet with $\GK(A)= \infty$.  Such examples do exist if one removes the graded or noetherian assumptions \cite[Remark after Example 1.4]{StZh97}.

\subsection{Auslander conditions}

The study of the interaction between homological conditions on an algebra and its ring-theoretic properties long predates the concept of an AS regular algebra.  In commutative algebra, one has such results as far back as the 1950's, such as the Auslander-Buchsbaum theorem that a regular local ring is a UFD.  The following is an important earlier notion  of regularity for noncommutative rings.
\begin{definition}
Let $A$ be a (not necessarily graded) noetherian
$\kk$-algebra $A$ with $\on{gl.dim}(A) = d < \infty$.
For an  $A$-module $M$, set $j(M) = \inf \{ i | \Ext^i_A(M, A) \neq 0 \}$, the \emph{$j$-number} or \emph{grade} of the module $M$.  $A$ is called \emph{Auslander regular} if for every finitely generated module $M$ and $0 \leq i \leq d$, all $A$-submodules $N$ of $\Ext^i_A(M,A)$ have $j(N) \geq i$.

If $\GK(A) < \infty$, the algebra $A$ is \emph{GK-Cohen-Macaulay} if for all finitely generated modules $M$, $\GK(M) + j(M) = \GK(A)$.
\end{definition}
At the time that AS regular algebras were introduced, the notion of Auslander regularity had already been important in ring theory for investigating examples related to rings of differential operators and Lie theory, for example.  Levasseur showed that if $A$ is a connected graded algebra with $\GKdim(A) < \infty$, then if $A$ is Auslander regular, it must be Artin-Schelter regular \cite[Theorem 6.3]{Lev92}.  Conversely, there is no known example of an AS regular algebra that is not also Auslander regular.

\subsection{Ring-theoretic properties of AS regular algebras}

Here are some basic questions about properties that AS regular algebras might satisfy.
\begin{question}
\label{ques:props}
Let $A$ be an AS regular algebra of global dimension $d$.
\begin{enumerate}
    \item Is $A$ noetherian? %If so, is $A$ strongly noetherian?
    \item Is $A$ a domain?  
    \item Does $\GK(A) = d$?     
    \item Is the minimal number of generators of $A$ as an algebra at most  $d$?
    \item Does $A$ have Hilbert series of the form $h_A(t) = \frac{1}{(1-t^{d_1}) \dots (1-t^{d_m})}$ for some $d_1, \dots d_m$ (the same as some graded commutative polynomial ring)? 
    \item Is $A$ Auslander regular?  If so, is $A$ GK-Cohen Macaulay?    
    \item Is $A$ a maximal order inside its Goldie quotient ring?
\end{enumerate}
\end{question}
The answer to each one of these questions is yes for every known example of regular algebra for which the property has been checked, and so the the answer is conjectured to be yes in every case.
Note that a maximal order is the noncommutative analog of an integrally closed ring in commutative algebra.

The following notion will appear occasionally below.  A $\kk$ algebra $R$ is a \emph{polynomial identity (PI) algebra} if there is some nonzero $f \in \kk \langle x_1, \dots, x_n \rangle$ such that $f(r_1, \dots, r_n) = 0$ for all $r_1, \dots, r_n \in R$.
If $R$ is commutative, then $R$ is a PI algebra by taking $f =x_1x_2 -x_2x_1$, and in general PI algebras have special properties that make them ``close to commutative."  For example, an algebra $R$ that is finitely generated as a module over its center $Z(R)$ is PI, but the converse is not true in general.

There are known connections between some of the questions above. Levasseur proved that if $A$ is connected graded and Auslander regular, then $A$ is a domain \cite[Theorem 4.8]{Lev92}.    In \cite[Section 3]{ATV2}, Artin, Tate and Van den Bergh proved that a regular algebra of dimension at most $4$ which is noetherian is a domain.    Stafford proved  that if an AS-regular algebra $A$ is noetherian, Auslander regular and GK-Cohen-Macaulay, then $A$ is indeed a maximal order.  Moreover, such an $A$ is a PI algebra if and only if it is a finite module over its noetherian center \cite[Theorem 2.10]{Sta94a}.  Stafford and Zhang examined the homological properties of graded PI algebras more generally in \cite{StaZh94b}.  They proved that if $A$ is a PI noetherian AS regular algebra, then $A$ is Auslander regular and GK-Cohen-Macaulay (so in particular, a domain and maximal order, as noted above).

One may also ask the opposite sort of question:  are there conditions weaker than those in the definition of AS regular that imply AS regularity, perhaps for special kinds of algebras?  
Stephenson and Zhang noted that connected graded noetherian algebras of global dimension $2$ are AS regular \cite[Theorem 0.4]{StZh97}.  In other words, the Gorenstein condition comes for free in the presence of the noetherian property (we already noted earlier that finite GK-dimension is forced).  They extended this result to some extent to dimension $3$, proving that any connected graded noetherian algebra $A$ of global dimension $3$ which is generated by 2 elements, has a graded free resolution of $\kk$ of the expected form, or has relations all in the same degree, must be AS regular \cite[Theorem 0.1]{StZh00}.  

Based on this evidence, one might conjecture that connected graded noetherian algebras of finite global dimension must be AS regular.  In \cite{RoSi12}, however, the author and Sierra found examples of algebras $A$ with $\gldim(A) = 4$ and $h_A(t) = \frac{1}{(1-t)^4}$ which are noetherian domains, but for which the Gorenstein condition fails.  These examples are even Koszul, as defined in Section~\ref{sec:extalgebra} below.

%These algebras have additional odd properties.  While they have $\GK(A) = 4$, and so projectively should correspond to a 3-fold, they behave more like a projective surface in the sense that the graded ring of quotients of $A$ is of the form $\kk(u,v)[z ; \sigma]$, i.e. geometrically they are ``birational" to a commutative surface in some sense.
%As we noted earlier, 
%noncommutative connected graded algebras with finite global dimension need not be noetherian or domains.  Algebras without polynomially bounded growth are ``big" in some sense, and the free algebra itself in more than one generator is not noetherian, so it is not surprising that some condition on growth is necessary if one expects to obtain noetherian algebras.  However, it is much more mysterious why the AS Gorenstein condition should force an algebra to be very well-behaved (for instance a domain).  Every known example of a regular algebra has a number of good properties, however, so it is natural to conjecture that these properties hold in general.  We review a number of these properties in this section and what little has been proved.

\subsection{Noncommutative Proj}

%If $A$ is a connected graded commutative algebra finitely generated in degree $1$, then $\Proj A$ is a projective scheme.  
When $A = \kk[x_1, \dots, x_n]$ is a polynomial ring, with some choice of grading given by assigning positive degrees (or \emph{weights}) to the variables, then the associated scheme $\Proj A$ is a called a \emph{weighted projective space}.  For this reason, intuitively a regular algebra $A$ can be thought of as a coordinate ring of a \emph{noncommutative} weighted projective space.  But what is the space?

Artin and Zhang showed that instead of a space, one can attach a category to $A$ which serves many of the same functions \cite{AZ1}.  We loosely describe the construction here.  Let $A$ be a noetherian connected graded finitely generated $\kk$-algebra.  Let $\rgr A$ be the abelian category of finitely generated $\mb{Z}$-graded right  $A$-modules, where the morphisms are graded (degree-preserving) module homomorphisms. 
%A finitely generated module $M \in \rgr$ is called \emph{torsion} in this context if for every $m \in M$ there is $n$ such that $ m A_{\geq n}= 0$; equivalently, if $M$ is finite-dimensional over $\kk$.  
Let $\rfd A$ be the full subcategory of $\rgr A$ consisting of modules which are finite dimensional over $\kk$.  Using a construction of quotient category due to Gabriel, one defines $\rproj A = \rgr A/\rfd A$ and calls this category the 
``noncommutative proj" associated to $A$.   If $A$ is commutative, generated in degree $1$, then the category $\rproj A$
is equivalent to the category of coherent sheaves on the usual scheme $X = \Proj A$.  So really we are generalizing the idea of a category of sheaves rather than a space to the noncommutative setting, but this turns out to be enough in order to develop a rich theory.  For example, the image of a point module $M$ for $A$ is a simple object in the category $\rproj A$, and so does behave geometrically as a ``point".

Noncommutative projective geometry is its own large subject.  See the survey article \cite{StVdB01}; the author's course notes \cite{Rog16} also give an introduction to the noncommutative proj construction.  The categories $\rproj A$ for regular algebras $A$ are very important in noncommutative geometry, because they represent noncommutative weighted projective spaces.  In particular, if $A$ is a quadratic regular algebra of dimension $3$, then $\rproj A$ is a noncommutative $\mb{P}^2$, and if $A$ is a cubic regular algebra of dimension $3$, then $\rproj A$ is a noncommutative $\mb{P}^1 \times \mb{P}^1$.  In this article, we will focus on the algebraic aspects of regular algebras, and will not attempt to survey the work that is focused primarily on their associated noncommutative projective spaces.  The interaction between the algebra $A$ and its category $\rproj A$ is an essential part of the story, however, and we will want to refer to the categories $\rproj A$ occasionally below.

\subsection{Twists}
\label{sec:twists}

Graded twists are a fundamental construction for graded rings that have played a role in much of the work on regular algebras, including the original work of Artin, Tate, and Van den Bergh.

Let $G$ be a semigroup.  Suppose that $A$ is a $G$-graded $\kk$-algebra, so $A = \bigoplus_{g \in G} A_g$ with $A_g A_h \subseteq A_{gh}$ for all $g, h \in G$.  Assume that for each $g \in G$ we have a $\kk$-linear bijection $\tau_g: A \to A$.  The collection  $\tau = \{ \tau_g | g \in G \}$ is called a \emph{twisting system} if $\tau_g(y \tau_h(z)) = \tau_g(y) \tau_{gh}(z)$ for all $g, h, l \in G$ and $y \in A_h$, $z \in A_l$.  Given a twisting system $\tau$ on $A$, one can define a new $G$-graded algebra $A^{\tau}$ which is called a \emph{twist} of $A$.  The algebra $A^{\tau}$ has the same underlying graded vector space as $A$, but a new multiplication $*$, defined on homogeneous elements $x \in A_g, y \in A_h$ by 
$x * y = x \, \tau_g(y)$, and extended linearly.  The condition defining a twisting system is exactly what is needed to ensure that $A^{\tau}$ is associative.  Zhang defined and studied this notion in \cite{Zha96}, and these twists are often called Zhang twists.  Zhang proved that for any twisting system $\{\tau_i \}$ on $A$, one has equivalent categories of graded modules $\rGr A \sim \rGr A^{\tau}$.
%and of course also equivalent categories of finitely generated graded modules $\rgr A \sim \rgr A^{\tau}$.

Consider in particular the case that $G = \mb{N}$, so $A$ is $\mb{N}$-graded.  In this case there is an easy way to construct 
twisting systems:  for any algebra automorphism $\sigma: A \to A$, if $\tau_n = \sigma^n$ for all $n$ then $\tau = \{ \tau_n \}$ is a twisting system.  The new multiplication in this case is given by $x * y = x \sigma^m(y)$ for all $x \in A_m, y \in A_n$.  As a converse to the result above in the $\mb{N}$-graded case,  if $A$ and $B$ are connected $\mb{N}$-graded and $A_1 \neq 0$, then $\rGr A \sim \rGr B$ implies that there exists a twisting system $\{ \tau_i \}$ on $A$ such that $B \cong A^{\tau}$ as graded algebras \cite[Theorems 3.1, 3.5]{Zha96}.  

Because the noncommutative projective scheme $\rproj A$ is a quotient category of the category of finitely generated graded $A$-modules, it is immediate that twisting does not change the underlying noncommutative scheme $\rproj A$.  However, in general twisting certainly does change the algebra up to isomorphism, so one can get a lot of different algebras as twists $A^{\tau}$, which can be thought of as different noncommutative coordinate rings of the same noncommutative scheme.  Zhang also proved that many properties are preserved under twist: in particular, if $A$ is AS regular, then so is $A^{\tau}$ for any twisting system $\tau = \{ \tau_i \}$ \cite[Theorem 5.11]{Zha96}.  Thus twisting gives an equivalence relation on the the set of AS regular algebras.  In classification problems for regular algebras, it is often useful to work up to twist equivalence.

A related kind of twist is a cocycle twist.  Suppose that $A$ is any algebra with a grading by a group $G$.  One may attempt to deform the multiplication of $A$ to a new one $*$ by choosing for each pair of elements $g, h \in G$ a nonzero constant $\mu(g, h) \in \kk$ such that for $a \in A_g, b \in A_h$, we have $a * b = \mu(g,h) ab$.  This gives a new associative multiplication on $A$ if and only if $\mu: G \times G \to \kk^{\times}$ is a $2$-cocycle.  The new algebra $(A, *)$ is written as $A_{\mu}$.  One may see that this is in fact a Zhang twist where each $\tau_g$ scales $A_h$ by $\mu(g, h)$, and in this case the 2-cocycle condition is equivalent to the twisting system axiom.

Now suppose that we have a connected graded algebra $A$ that has an additional grading by a group $G$ which refines the given $\mb{N}$-grading.  Then given a $2$-cocycle $\mu$ for $G$, the new algebra $A_{\mu}$ will still be $\mb{N}$-graded.  However, the preservation of nice properties under this kind of twist no longer follows from Zhang's results.  Davies proved that AS regularity and other good properties are in fact invariant under twisting when $G$ is a finite abelian group \cite[Theorem 1.2]{Dav16}.  He also showed that some families of examples of $4$-dimensional regular algebras classified by the author and Zhang in \cite{RoZh12} are related by non-obvious twists of this kind.  
%In any classification of regular algebras one would ideally also want to understand the equivalence classes of cocycle twists.

%\textcolor{blue}{known results on twists?}

\section{Algebras related to an elliptic curve}

The Sklyanin algebras of dimension $3$, defined above in Example~\ref{ex:sklyanin3}, are the most interesting family of 3-dimensional regular algebras, due to the nontrivial and beautiful connections between these algebras and the geometry of elliptic curves.  In fact, they belong to a large family of examples related to elliptic curves, as we describe in this section.

\subsection{The Sklyanin algebras of dimension $4$}

The first regular algebras of dimension $4$ to be studied in detail were the $4$-dimensional Sklyanin algebras.  
%One can write a down-to-earth presentation of this family of algebras:
\begin{example}
\label{ex:skl4}
Let $\alpha, \beta, \gamma \in \kk$ with $\alpha + \beta + \gamma + \alpha \beta \gamma = 0$.

The $4$-dimensional Sklyanin algebra is defined by 
\[
S(\alpha, \beta, \gamma) = \kk \langle x_0, x_1, x_2, x_3 \rangle/(r_1, \dots, r_6)
\]
where the six relations are
\begin{gather*}
r_1 = x_0x_1 -x_1x_0 - \alpha(x_2x_3 + x_3x_2), \qquad r_2 = x_0x_1 + x_1x_0 - (x_2x_3 - x_3x_2), \\
r_3 = x_0x_2 - x_2x_0 - \beta(x_3x_1 + x_1x_3), \qquad r_4 = x_0x_2 + x_2x_0 - (x_3x_1 - x_1x_3), \\
r_5 = x_0x_3 - x_3x_0 - \gamma(x_1x_2 + x_2x_1), \qquad r_6 = x_0x_3 + x_3x_0 - (x_1x_2 - x_2x_1).
\end{gather*}
\end{example}

The algebra $S(\alpha, \beta, \gamma)$ was proved to be AS regular precisely when $\{\alpha, \beta, \gamma \}$ does not contain the subset $\{1, -1\}$, by Smith and Stafford \cite{SmSt92}.  They proved that such a regular $S$ is a noetherian domain with Hilbert series $H_S(t) = \frac{1}{(1-t)^4}$.  The proof uses a similar idea to the proof that the Sklyanin algebras of dimension $3$ are regular in \cite{ATV1}.  In more detail, $S$ has two linearly independent central elements of degree two, $g_1$ and $g_2$, say, such that $A/(g_1, g_2) \cong B = B(E, \mc{L}, \sigma)$ is a twisted homogeneous coordinate ring on an elliptic curve $E$.  The properties of $B$ can be determined geometrically, and good properties can be lifted from $B$ to $S$.   In this case $E$ is not the entire point scheme $X$, which also has $4$ isolated points in addition to $E$.  Further properties of $S$ were proved by Staniszkis \cite{Stan94} and Levasseur and Smith \cite{LeSm93}.  

The $4$-dimensional Sklyanin algebras $S$ have also been a good starting place for constructing other interesting examples.  Since $S$ is presented by $4$ generators and $6$ quadratic relations, $B$ is presented by $4$ generators and $8$ quadratic relations.  Stafford examined another class of rings obtained by removing a different choice of $2$ relations from the presentation of $B$, and showed these are also AS regular of dimension $4$ \cite{Sta94b}. This showed that the ring $B$ does not uniquely determine the $4$-dimensional regular algebra mapping onto it, unlike the situation for regular algebras of dimension $3$.

In \cite{Dav17}, Davies constructed regular cocycle twists of $4$-dimensional Sklyanin algebras $S$ that have quite different properties from $S$.  For example, these twists have only finitely many point modules up to isomorphism, unlike $S$ which has an entire elliptic curve of point modules, demonstrating how a cocycle twist can change the geometry associated to an algebra dramatically, unlike a usual Zhang twist.  Chirvasitu and Smith investigated many additional properties of these and more general kinds of twists in a series of papers \cite{ChSm19} \cite{ChSm17}.

%Since $S$ has a $2$-dimensional space of central elements, it is also natural to investigate the algebras $S/(\Omega)$ given by modding out by one central element $0 \neq \Omega \in \kk g_1 + \kk g_2$, since these algebras are intermediate between $S$ and $B$.  These are AS Gorenstein rather than regular, and can be thought of as coordinate rings of noncommutative quadric surfaces.  In \cite{VdB96}, Van den Bergh proved a translation principle for these which gives equivalences between certain of the corresponding noncommutative projective schemes $\rproj S/(\Omega)$.  This translation principle is an analog of similar results about the primitive quotient rings of the enveloping algebra $U(\mf{sl}_2)$.  
%See also \cite{SmVdB13} for further results about noncommutative quadric surfaces.

\subsection{Feigin and Odesskii's Elliptic Algebras}
\label{sec:elliptic}

Roughly at the same time that Artin, Schelter, Tate, and Van den Bergh were laying down the foundations for the classification of regular algebras of dimension $3$, Feigin and Odesskii defined a class of algebras built out of the geometry of an elliptic curve \cite{OdFe89}.  Serendipitously, these algebras turned out to be important examples of AS regular algebras.

There are several minor variations on the basic definition that have appeared in the literature; we follow here the presentation given by Chirvasitu, Kanda, and Smith in \cite{CKS21a}.

Let $E$ be a complex elliptic curve, defined as $\mb{C}/\Lambda$ for a lattice $\Lambda = \mb{Z} + \mb{Z} \eta$, where $\eta$ is in the upper half of the complex plane.  Fix $1 \leq k < n$ with $\gcd(k, n) = 1$ and a point $\tau \in \mb{C} - (1/n)\Lambda$.  A \emph{theta function} is a holomorphic function which is quasi-periodic with respect to the lattice, in the sense that translating by one of the lattice generators $1$ or $\tau$ scales the function by a specified complex exponential function.  More specifically, here we are interested in holomorphic functions $f: \mb{C} \to \mb{C}$ which satisfy 
\[
f(z+1) = f(z)\ \qquad \text{and}\ \qquad f(z + \eta) = -e^{-2\pi i nz}f(z). 
\]
These are called \emph{theta functions of order $n$} as they have $n$ zeros inside every period parallelogram of the lattice.  The set of all such theta functions of order $n$ is a vector space over $\mb{C}$ of dimension $n$ denoted $\Theta_{n}(\Lambda)$, which is an irreducible representation of the finite Heisenberg group $H_n$ of order $n^3$ in a natural way. One chooses a basis $\theta_0, \ldots, \theta_n$  of $\Theta_n(\Lambda)$ which is well-suited to the action of $H_n$; see \cite[Proposition 2.6]{CKS21a}.  In the next definition we take the indices of the generators $x_i$ and the theta functions $\theta_i$ modulo $n$.

\begin{definition}
Define
\[
Q_{n,k}(E, \tau) = \kk \langle x_0, x_1, \dots, x_{n-1} \rangle/( r_{ij} | 1 \leq i,j \leq n)
\]

where
\[
r_{ij} = \sum_{r \in \mb{Z}_n} \frac{\theta_{j-i + (k-1)r}(0)}{\theta_{j-i-r}(-\tau)\theta_{kr}(\tau)} x_{j-r} \, x_{i+r},
\]
for $\theta_0, \dots \theta_{n-1}$ the special basis of $\Theta_n(\Lambda)$ chosen above.
\end{definition}

Clearly there is a lot to unpack in this definition, and in a survey article like this we cannot give much of a feel for these algebras.  The already mentioned paper of Chirvasitu, Kanda, and Smith \cite{CKS21a} is the recommended starting point, as it is the first paper in a large project by these authors to provide more explicit details of the basic properties of these algebras, many originally stated by Feigin and Odesskii with brief or no proofs.  While we assumed $\tau \not \in (1/n)\Lambda$ above, there is a natural way to extend the definition so that an algebra $Q_{n,k}(E, \tau)$ is defined for every $\tau \in \mb{C}$.

Feigin and Odesskii were motivated by algebras defined by Sklyanin in papers predating any of the work on regular algebras \cite{Skl82}, \cite{Skl83}, which studied certain algebras that were defined in terms of Baxter's ``elliptic solution" to the quantum Yang-Baxter equation.  In the more general setup of Feigin and Odesskii, the algebras defined by Sklyanin are the algebras $Q_{4,1}(E, \tau)$, and of course this is why the general construction is named after Sklyanin.  For generic $\tau$, these are Sklyanin algebras of dimension $4$ whose presentation was given in Example~\ref{ex:skl4}, in the case $\kk = \mb{C}$.  Similarly, for generic $\tau$ the algebras $Q_{3,1}(E, \tau)$ are the 3-dimensional Sklyanin algebras in Example~\ref{ex:sklyanin3}.  Thus the algebras $Q_{n,1}(E, \tau)$ are direct generalizations of the algebras in Sklyanin's papers, and are commonly referred to as Sklyanin algebras of dimension $n$.

In related work, Tate and Van den Bergh gave an essentially different approach  to the algebras $Q_{n,1}(E, \tau)$ from Feigin and Odesskii's \cite{TVdB96}.  In particular, they showed these algebras are noetherian domains.  See also \cite{Teo96}.  Smith studied the point modules of these algebras in \cite{Smi94b}.  More generally, all linear modules of $Q_{n,1}(E, \tau)$ were classified by Staniszkis \cite{Stan96}.  As in the case of $n = 3$ that appeared in \cite{ATV1}, considering point modules for $Q_{n,1}(E, \tau)$ gives a natural embedding of $E$ inside $\mb{P}^{n-1}$.  Staniszkis showed that the linear modules of $Q_{n,1}(E, \tau)$ correspond to those linear subspaces of $\mb{P}^{n-1}$ which relate to the embedded elliptic curve $E$ in a certain way.  She proved that these linear modules are GK-Cohen-Macaulay and critical with respect to GK-dimension.  Zhang \cite{Zha98b} showed that the $d$-linear modules of $Q_{n,1}(E, \tau)$ have projective dimension $n-d$ and have linear free resolutions.  

As we see above, the Sklyanin algebras---that is, the $Q_{n,k}(E, \tau)$ with $k = 1$---have been the most studied among Feigin and Odesskii's algebras.  As already mentioned, Chirvasitu, Kanda, and Smith have a recent series of papers  \cite{CKS-arxiv19} \cite{CKS21a} \cite{CKS21b}   \cite{CKS-arxiv21} \cite{CKS23} revisiting the work of Feigin and Odesskii on the algebras $Q_{n,k}(E, \tau)$, with more attention on the case $k \neq 1$.  In more detail, after an introductory paper \cite{CKS21a} they define and study the characteristic variety $X_{n/k}$, which is a quotient of $E^g$ by a finite group---here $g$ is a positive integer depending on a continued fraction expansion of $n/k$ \cite{CKS-arxiv19}.  When $k=1$, $X_{n/k} = E$.  There is a map from $Q_{n,k}(E, \tau)$ to a certain canonical twisted homogeneous coordinate ring $B(X_{n/k}, \mc{L}, \sigma)$ which is surjective in many cases \cite{CKS21b}.  In particular, they verify carefully that $Q_{n,k}(E, \tau)$ is a regular algebra with Hilbert series $\frac{1}{(1-t)^n}$ for all but countably many $\tau$ \cite{CKS23}.  There are still many basic open questions in this area, and we close this section with just a few.  Are the regular algebras of the form $Q_{n,k}(E, \tau)$ noetherian domains when $k \neq 1$?  What are the finite-dimensional representations over $Q_{n,k}(E, \tau)$?  Can one define analogs of the algebras $Q_{n,k}(E, \tau)$ for any abelian variety $E$?

\section{Linear modules over regular algebras}

Point modules, and the scheme parametrizing them, were crucial to the geometric classification of AS regular algebras of dimension $3$ given in \cite{ATV1}.  Naturally, many investigations of regular algebras of dimension $4$ have also considered the point scheme, with a hope that one might be able to recover the algebra from geometric data on the point scheme in some cases, for example.  Starting in dimension $4$, larger dimensional linear modules such as line modules are also natural objects of study.  

\subsection{Hilbert schemes}
\label{sec:hilb}

We first review the technical theory underpinning the generalization of the point scheme to parameter spaces of higher-dimensional modules.  As we saw earlier, for any connected graded algebra $A$ generated in degree $1$, the point modules are in bijection with an inverse limit of schemes $\invlim \, X_d$, where $X_d$ is a projective scheme parametrizing truncated point modules of length $d+1$.  This construction is most useful when this inverse limit stabilizes, in which case there is a projective scheme parametrizing the point modules.

More generally, one would like to find a nice moduli space for the isomorphism classes of cyclic graded modules with a given fixed Hilbert series---a \emph{Hilbert scheme}.  Schemes which parametrize point modules, line modules, and other linear modules are then just a special case.   Artin and Zhang developed in \cite{AZ2} such a theory of noncommutative Hilbert schemes.   They showed that if a connected graded algebra $A$, finitely generated in degree $1$, is \emph{strongly noetherian}, that is, $A \otimes_{\kk} R$ is noetherian for all commutative noetherian $\kk$-algebras $R$, then the Hilbert schemes are projective schemes.  In particular, in this case there is a projective point scheme, line scheme, and so on.  Commutative noetherian finitely generated graded algebras are automatically strongly noetherian, and Artin and Zhang's result simply recovers the commutative theory of Hilbert schemes in that case.  However, noncommutative noetherian connected graded algebras may fail to be strongly noetherian, and the inverse limit parametrizing their point modules may fail to stabilize \cite{Jor01} \cite{Rog04}.  A large class of such examples called \emph{na\"ive blowups} can be constructed as certain subalgebras of twisted homogeneous coordinate rings \cite{KRS05}.

All known regular algebras are strongly noetherian.  In particular examples, though, one tends to show that the point modules or line modules are parametrized by a projective scheme directly, without first attempting to prove the algebra is strongly noetherian.  In fact, in regular algebras the inverse limits which parametrize point and line modules tend to stabilize in the lowest possible degree, similar to the behavior seen in regular algebras of dimension $3$.  For simplicity, let us focus on the case of a regular algebra $A$ of dimension $d$ with quadratic relations and Hilbert series $\frac{1}{(1-t)^d}$.  In this case, in all known examples the inverse limit $\invlim \, X_d$ which parametrizes point modules stabilizes already at $X_2 \subseteq \mb{P}^{m-1} \times \mb{P}^{m-1}$, and moreover, $X_2$ is the graph of an automorphism $\sigma: X \to X$ for some $X \subseteq \mb{P}^{m-1}$. When this happens, $X \cong X_2$ parametrizes the point modules, and truncation shift is an automorphism on the set of point modules (given by $\sigma$).  

Shelton and Vancliff gave an extensive study of line modules for AS regular algebras $A$ of dimension $4$ with quadratic relations.  They showed that if $A$ is Auslander regular and GK-Cohen Macaulay, then the line modules of $A$ are parameterized by a projective scheme that can easily calculated as a subscheme of a Grassmannian \cite[Lemma 2.4, Remark 2.10]{ShVa02a}.  Under the same condition, they showed that this line scheme has dimension at least $1$, and that every point module is a quotient of some line module \cite[Proposition 3.1]{ShVa02a}.  They computed the line scheme explicitly in a number of interesting examples \cite{ShVa02b}.  In more recent work, Chirvasitu, Smith, and Vancliff studied the geometric properties of the scheme parametrizing line modules when it has dimension $1$, proving that it always is a subscheme of degree $20$ inside an ambient $\mb{P}^5$ \cite{CSV20}.

\subsection{Questions about the point scheme}

Here are some questions related to the point modules of regular algebras.  Again, for convenience we focus on the case of algebras with quadratic relations, but one can ask similar questions in general.
\begin{question}
\label{ques:props2}
Let $A$ be an AS regular algebra, generated in degree $1$ with $\dim_{\kk} A_1 = m$.
\begin{enumerate}
    \item Assuming $A$ is noetherian, is $A$ strongly noetherian? 
    \item Suppose that $A$ has quadratic relations, and let $X_2 \subseteq \mb{P}^{m-1} \times \mb{P}^{m-1}$ be the subscheme defined by the multilinearizations of the relations.  Does the inverse limit $\invlim \, X_d$ parametrizing the point modules already stabilize at $X_2$?  
    If so, is $X_2$ the graph of an automorphism of some subscheme of $\mb{P}^{m-1}$?
    \item When is the algebra $A$ determined by its point scheme data?  If it is not, what if we include the data of some other Hilbert schemes?
    \item Is whether or not $A$ is a polynomial identity algebra (PI) related in a nice way to the point scheme data?  If it is not, what if we include the data of some other Hilbert schemes?
    \end{enumerate}
\end{question}
These questions are inspired by what happens in dimension $3$.  For any AS regular algebra $A$ of dimension $3$, the answer to (1) is yes by a similar argument as the proof that $A$ is noetherian.  For a quadratic AS regular algebra $A$ of dimension $3$, the scheme $X_2$ is indeed the graph of an automorphism $\sigma: E \to E$ as in (2), and for (3) the data $(E, \mc{O}_{\mb{P}^2} \vert_E, \sigma)$ determines the algebra, as we saw earlier  in Theorem~\ref{thm:ATV-converse}.  For (4), \cite[Theorem 7.1]{ATV2} proved that $A$ is PI if and only if the automorphism $\sigma$ of the point scheme has finite order.   

The answer is also yes to questions (1) and (2) in all known higher-dimensional examples for which it has been checked.  Assume in the following discussion that $A$ is regular of dimension $4$, generated in degree $1$ with quadratic relations and Hilbert series $\frac{1}{(1-t)^4}$.  It is known that if $A$ is Auslander regular and GK-Cohen Macaulay, then (2) holds; this was first proved by Vancliff and Van Rompay in \cite[Theorem 1.10]{VaVR97}, under an additional minor hypothesis which was later removed by Shelton and Vancliff \cite[Theorem 1.4]{ShVa99b}.  If (2) holds, one may consider the subspace of $A_1 \otimes_{\kk} A_1$ consisting of elements whose multilinearizations vanish along $X_2 \subseteq \mb{P}^{m-1} \times \mb{P}^{m-1}$.  In nice cases, this may be equal to the span of the multilinearlizations of the degree $2$ relations of $A$, and hence $A$ is determined by $X_2$, giving an answer to (3).  Similarly, it is possible for $A$ to be completely determined by the line scheme  data of $A$; in particular, this happens if the line scheme is of dimension $1$, which is known to be the minimal possible dimension for the line scheme \cite[Theorem 4.1]{ShVa02a}.

%Following \cite[Definition 2.4]{Mor06a} we say that $A$ is \emph{geometric} in this case.

\subsection{The dimension of the point scheme of regular algebras of dimension $4$}

Suppose that $A$ is a regular algebra of dimension $4$ with quadratic relations and Hilbert series 
$\frac{1}{(1-t)^4}$.  In general, the larger the point scheme $X$ of the regular algebra $A$, the more information it is expected to provide.  At one extreme, it could be that $X \cong \mb{P}^3$.  This case is rather trivial, as $A$ will then be a twisted homogeneous coordinate ring $B(\mb{P}^3, O(1), \sigma)$, which is just a Zhang twist of a polynomial ring in $4$ variables.  

The next is case where $\dim X = 2$.
One might expect this case to bear some similarity to the most typical situation for quadratic regular algebras of dimension $3$, in which the point scheme is a degree three divisor inside $\mb{P}^2$.
A number of cases with $\dim X = 2$ have been studied in detail.  In \cite{Van94}, Vancliff defined a natural class of such algebras $A$ with a normal element $x \in A_2$ such that $A/xA \cong B(Q, \mc{L}, \sigma)$, for some nonsingular quadric hypersurface $Q$ in $\mb{P}^3$ with automorphism $\sigma$ and invertible sheaf $\mc{L} = \mc{O}(1) \vert_{Q}$.  The quadric $Q$ is not quite the point scheme $X$ of $A$, but rather $X = Q \cup L$ for some line $L$ intersecting $Q$ in two points.  Vancliff wrote down the relations of all such algebras; they include interesting examples such as the coordinate ring of quantum $2 \times 2$ matrices $\mc{O}_q(M_2(\kk))$ for $q^2 \neq 1$.  She also examined the line modules and their incidence relations with the point modules.  

Van Rompay \cite{VRo96} constructed similar examples using Segre products.  Then in \cite{VaVR97}, Vancliff and Van Rompay classified all regular algebras $A$ of dimension $4$ which surject onto a twisted homogeneous coordinate ring $B(Q, \mc{L}, \sigma)$ for a nonsingular quadric $Q$; in many cases the point scheme of $A$ is of the form $X = Q \cup L$, where in general the line $L$ can be tangent to $Q$ or embedded in $Q$ as a line of multiplicity $2$.  They showed that it is also possible to have $X = X_2 = Q$, and in this case there is a whole family of regular algebras that surject onto $B(Q, \mc{L}, \sigma)$; that is, the point scheme data does not determine the relations of the algebra \cite{VaVR97}. However, in \cite{VaVR00}, they give a further analysis of this class of algebras using line modules, and prove that the schemes of point and line modules together determine the relations of the algebra, giving a kind of positive answer to Question~\ref{ques:props2}(3) in this case.  

Shelton and Vancliff \cite{ShVa99b} extended this program to singular quadrics $Q$ of rank 3 in $\mb{P}^3$, classifying all regular algebras $A$ of dimension $4$ which surject onto $B(Q, \mc{L}, \sigma)$ for such $Q$.

%Shelton and Vancliff also show that for some examples of regular algebras of dimension $4$, for which the point scheme $X_2$ does not determine the relations of the algebra, the point scheme together with the line scheme does determine the algebra, which gives a kind of answer to question (3) above in a special case \cite[Example 4.5]{ShVa02a}.

We are unaware of any systematic study of $4$-dimensional regular algebras with point scheme $X$ such that $\dim X = 1$, though of course there are many examples, notably the $4$-dimensional Sklyanin algebras.

The final case is when the point scheme $X$ has $\dim X = 0$, so that $X$ is a collection of points, possibly with non-reduced scheme structure.  
In contrast to the examples where the point scheme is a quadric surface, Shelton and Vancliff proved the surprising result that in fact the point scheme $X$ data determines the relations of a regular algebra $A$ of dimension $4$ as above, whenever $\dim X = 0$ \cite[Theorem 4.1, Remark 4.2]{ShVa02a}.  So Question~\ref{ques:props2}(3) has a nice answer in this case.   There is a useful construction, which we review next, that produces many examples of this type with $X$ of various different configurations.

\subsection{(Skew) graded Clifford algebras}
\label{sec:clifford}

Fix $n \geq 2$.  Fix nonzero scalars $\mu_{ij} \in \kk$ such that $\mu_{ii} = 1$ for all $i$ and 
$\mu_{ij}\mu_{ji} = 1$ for all $i \neq j$.  A matrix $M \in M_n(\kk)$ is $\emph{$\mu$-symmetric}$ if 
$M_{ij} = \mu_{ij} M_{ji}$ for all $i,j$.

Suppose we are given a collection of $n \times n$ $\mu$-symmetric matrices $M_1, \dots, M_n$ in $M_n(\kk)$.  We define an algebra $A(\mu, M_1, \dots, M_n)$ by the following presentation:
\[
\kk \langle x_1, \dots, x_n, y_1, \dots, y_n \rangle /( x_ix_j + \mu_{ij} x_j x_i - \sum_{k=1}^n (M_k)_{ij} \, y_k \, | \, 1 \leq i, j \leq n ).
\]
Here the generators $x_i$ are of degree $1$ and the $y_i$ are of degree $2$; so the algebra is naturally connected $\mb{N}$-graded.  Recall that an element $x$ in a ring $R$ is \emph{normal} if $xR = Rx$, so that the right ideal $xR$ is an ideal.  Following Cassidy and Vancliff, $A = A(\mu, M_1, \dots, M_n)$ is called a \emph{graded skew Clifford algebra} if in addition there is a $\kk$-basis of $\kk y_1 + \dots + \kk y_n$, say $r_1, \dots, r_n$, which is a normalizing sequence in $A$ in the sense that each $r_i$ is normal in the factor ring $A/(r_1, \dots, r_{i-1})$.  While by definition these algebras have generators in degrees $1$ and $2$, the point of this construction is to find interesting regular algebras that are generated in degree $1$.  As long as the matrices $M_1, \dots, M_n$ are linearly independent over $\kk$, one may easily check that the generators $y_i$ may be eliminated from the presentation, and so $A$ is generated by $x_1, \dots x_n$.

If $\mu_{ij} = 1$ for all $i, j$, then the matrices $M_i$ are symmetric.  In this case we omit $\mu$ from the notation,  writing $A = A(M_1, \dots, M_n)$, and call this a \emph{graded Clifford algebra} if in addition $y_1, \dots, y_n$ are central in $A$ for all $i$.  This was the first case to be studied.
In this case, each matrix $M_k$ determines a quadric hypersurface $Q_k$ in $\mb{P}^{n-1} = \Proj \kk[z_1, \dots, z_n]$ given by the vanishing of $\sum_{i,j} (M_k)_{ij} z_i z_j = 0$.  The system of quadrics $Q_1, \dots, Q_n$ is called \emph{basepoint free} if $Q_1 \cap \dots \cap Q_n = \emptyset$.  Le Bruyn \cite{LeB95} proved that $A(M_1, \dots, M_n)$ is a regular algebra with Hilbert series $\frac{1}{(1-t)^n}$ if and only if the associated system of quadrics is basepoint free.  

\begin{example}[Section 2, \cite{VVW98}]
\label{ex:20pts}
Let $M_1, M_2, M_3, M_4 \in M_n(\kk)$ be the symmetric matrices such that 
\[ M(\alpha_1, \alpha_2, \alpha_3, \alpha_4) = \sum_{i=1}^4 \alpha_i M_i = \begin{bmatrix} 2\alpha_1 & \alpha_3  & \alpha_4 & \alpha_2 \\
\alpha_3 & 2 \alpha_2 & \alpha_1 & \alpha_3 + \alpha_4 \\
 \alpha_4 & \alpha_1 & 2 \alpha_3 & \alpha_1 \\
\alpha_2  & \alpha_3 + \alpha_4 & \alpha_1 & 2 \alpha_4 
\end{bmatrix}
\]
for indeterminates $\alpha_i$.  In this case the relations of the Clifford algebra $A = A(M_1, M_2, M_3, M_4)$ include $y_i = x_i^2$ for all $i$, so eliminating the generators $y_i$ we see that $A$ has the following presentation:
\begin{gather*}
A = \kk \langle x_1, x_2, x_3, x_4 \rangle/(r_1, r_2, r_3, r_4, r_5, r_6) \ \text{where} \\
r_1 = x_1x_4 + x_4x_1 - x_2^2, \qquad r_2 = x_2x_4 + x_4x_2 -x_3^2-x_4^2, \\
r_3 = x_3x_4 + x_4x_3 - x_1^2, \qquad r_4 = x_1x_3 + x_3x_1 - x_4^2, \\
r_5 = x_2x_3 + x_3x_2 - x_1^2, \qquad r_6 = x_1x_2 + x_2x_1 - x_3^2.
\end{gather*}
The basepoint-free condition on the associated quadratic forms is easily checked, so the Clifford algebra $A$ is regular.
\end{example}

Cassidy and Vancliff defined and studied the more general notion of graded skew Clifford algebra in \cite{CaVa10}.  Replacing $\kk[z_1, \dots, z_n]$ by the skew polynomial ring $S = \kk \langle z_1, \dots, z_n \rangle/(z_iz_j + \mu_{ij} z_j z_i | 1 \leq i < j \leq n)$, they defined an analog of the system of quadrics associated to $M_1, \dots, M_n$, where 
the quadric $Q_k$ associated to $M_k$ is now a hypersurface inside the point scheme of $S$.  They proved that the graded skew Clifford algebra $A(\mu, M_1, \dots, M_n)$ is AS regular if and only if the associated system of quadrics is base point free in a suitable sense, generalizing Le Bruyn's result.

As evidence that the skew graded Clifford algebra construction is quite general, Nafari, Vancliff, and Jun Zhang revisited the classification of regular algebras of dimension $3$, and showed that the majority of these algebras are, up to twist equivalence, either skew graded Clifford algebras of dimension $3$, or Ore extensions of skew graded Clifford algebras of dimension $2$ \cite{NVZ11}.
The reader can find a survey with more details on skew graded Clifford algebras in \cite{Vee17}.

\subsection{$0$-dimensional point schemes}
Consider a regular algebra $A$ of dimension $4$ with quadratic relations and Hilbert series $\frac{1}{(1-t)^4}$.  Van den Bergh observed that when the point scheme $X$ of $A$ is $0$-dimensional, then $X$ consists of $20$ points, counted with multiplicity \cite{VdB88}.  An outline of the argument can be found in \cite[p. 377]{Van15}.  See also \cite{ChKa23}, where Chirvasitu and Kanda generalize the argument to apply to regular algebras of dimension $4$ with other Hilbert series.

Van den Bergh's result raised the question of what kinds of finite schemes of length 20 actually occur for regular algebras with a $0$-dimensional point scheme.  For a regular graded Clifford algebra $A(M_1, \dots, M_4)$ as described above, the number $m$ of distinct closed points in the point scheme $X$ of $A$  can be efficiently calculated from the properties of the matrices $M_i$ \cite[Theorem 1.7]{VVW98}.  Namely, $m = 2r_2 + r_1$, where $r_i$ is the number of matrices of rank $i$ in the projectivized linear system $\mb{P}(\kk M_1 + \dots + \kk M_4)$ spanned by the defining matrices inside $M_4(\kk)$.   
For example, in Example~\ref{ex:20pts} this amounts to finding the points $(\alpha_1 : \alpha_2 : \alpha_3 : \alpha_4) \in \mb{P}^3$ where the matrix $M(\alpha_1, \alpha_2, \alpha_3, \alpha_4)$ has rank $\leq 2$. in \cite[Section 2]{VVW98} the authors show there are 20 distinct points in this example.  They also give an example which has one single point of multiplicity $20$.  

To calculate the values where a symmetric matrix with indeterminate entries has rank $\leq 2$ generally involves intersecting 10 cubic equations given by the vanishing of all of the $3 \times 3$ minors. Stephenson and Vancliff showed that in the setting of regular graded Clifford algebras, one can replace this calculation with the intersection of just $2$ cubic divisors \cite{StVa07}.  Using this idea, 
they fully analyzed the possible numbers $m$ of distinct points occurring for graded Clifford algebras with finite point scheme, showing that all $1 \leq m \leq 20$ are possible except for $m = 2, 15, 17, 19$.  

Much of this story generalizes nicely to skew graded Clifford algebras. Cassidy and Vancliff found a family of skew Clifford algebras $A(\mu, M_1, M_2, M_3, M_4)$, which as $\mu$ only is varied, produces examples of skew graded Clifford algebras with various different numbers of distinct points \cite[Table 1]{CaVa10}, demonstrating the considerable additional generality of the skew case.  In particular, they give examples in this way with $15$ or $17$ distinct point modules.  Vancliff and Veerapen gave a general formula for the number of distinct points in terms of the defining data for graded skew Clifford algebras, analogous to the formula for graded Clifford algebras \cite{VaVe14}.

\subsection{Other results on points}

In addition to point modules, so-called ``fat" point modules of regular algebras $A$ can also give interesting information.  A graded module $M$ over a connected graded algebra $A$ is called a \emph{fat point module} of multiplicity $d$ if $M$ is generated in degree $0$, has $\dim_{\kk} M_n = d$ for all $n \geq 0$, and every proper factor module of $M$ is finite-dimensional.  Just like point modules, fat point modules correspond to simple objects in the category $\rproj A$.  The moduli theory of fat point modules was studied by Daniel Chan \cite{Cha12}.  

Consider Question~\ref{ques:props2}(4). It was certainly not expected that the PI property for a regular algebra $A$ of dimension $4$ would always be determined just by its point scheme data.  Indeed, Stephenson and Vancliff gave examples with $0$-dimensional point scheme whose associated truncation shift automorphism has finite order, yet $A$ is not finite over its center, in contrast to what happens in dimension $3$  \cite{StVa06}.  Goetz and Shelton \cite{GoeSh06} found fat point modules of multiplicity $2$ for some of these regular algebras, and showed that truncation shift has infinite order on the set of fat points of multiplicity $2$.  This suggests that question (4) might still have some reasonable answer for regular algebras of dimension $4$ if fat point schemes are also included.  Goetz studied fat points for some other regular algebras of dimension $4$ as well \cite{Goe07}.

For regular algebras of dimension higher than $4$, there are few general results about point or line modules. 
Already in dimension $5$, something rather different happens with the point scheme.   Recall that for regular algebras of dimension $4$ with quadratic relations, the minimal dimension of the point scheme $X$ is $\dim X = 0$, and the point scheme data contains enough structure in this case to recover the relations of the algebra \cite[Theorem 4.1, Remark 4.2]{ShVa02a}.  It has long been known, on the other hand, that the same arguments used by Van den Bergh in \cite{VdB88} show that a generic quadratic regular algebra of dimension $5$ has a point scheme which is \emph{empty}.   The referee of this article noticed that an explicit example of a dimension $5$ regular algebra with empty point scheme was missing from the literature.  Vancliff has provided the construction of a regular Clifford algebra with such a property in a companion article for this volume \cite{Van23}.

Of course, point modules may still be interesting for special kinds of higher-dimensional regular algebras.  Belmans, De Laet, and Le Bruyn showed \cite{BDL16} that the point scheme for a quantum polynomial algebra $A = \kk \langle x_1, \dots, x_n \rangle/(x_j x_i - q_{ij} x_i x_j)$ has an interesting combinatorial structure, for arbitrary dimension $n$.  They prove it is a union of linear subspaces of $\mb{P}^{n-1}$ and determine its irreducible components.  Not all unions of linear subspaces can occur as such a point scheme; they give an additional necessary condition which is sufficient for small $n$.

\section{Extensions}

As we mentioned already in the introduction, an easy way to produce examples of regular algebras 
is by taking iterated graded Ore extensions.  This construction does not capture many interesting examples of regular algebras, however.  In this section we describe some more general constructions which have wider scope.  

\begin{comment}

As we have already mentioned above, if $A$ is a connected graded ring, then an Ore extension $A' = A[x; \sigma, \delta]$ will also be graded, with $\deg(x) = 1$, as long as $\sigma$ has degree $0$ and $\delta$ has degree $1$.

In this case it is known that $A$ is AS regular if and only if $A'$ is.
Moreover, if we start with an algebra $A$ which has other good properties, such as being strongly noetherian, Auslander regular and Cohen-Macaulay, then $A'$ inherits all of these properties \cite[Lemma 5.3]{ZhZh09}.  
In this way we can start with the regular algebra $\kk$ of dimension $0$ and obtain many examples of iterated Ore extensions 
\[
\kk[x_1][x_2; \sigma_2, \delta_2] \dots [x_m; \sigma_m, \delta_m]
\]
which are AS regular, and have all the good properties we expect of regular algebras, as in Questions~\ref{ques:props}.   
These are considered the easy examples of regular algebras in some sense, while examples such as the Sklyanin algebras which do not arise in this way are of greater interest.    That is not to say that we can easily classify all possible iterated Ore extensions, i.e. say explicitly what all of the possible choices of $\sigma_i$ and $\delta_i$ are, for large $m$.  But once we know that a regular algebra is an iterated Ore extension, further analysis of its properties becomes easier.
    
\end{comment}

\subsection{Normal extensions}

%A sequence $x_1, \dots, x_n \in A$ is a \emph{normalizing sequence} if the image of $x_i$ is normal in the algebra $A/(x_1, \dots, x_{i-1})$ for each $i \geq 1$.  
We say that a ring $B$ is a \emph{normal extension} of $A$ if there is a normal non-zero-divisor  $x \in B$ such that $B/xB \cong A$.   If  in fact $x$ is central in $B$, we call $B$ a \emph{central extension} of $A$.  
%We are only interested here in the case that $A$ and $B$ are connected $\mb{N}$-graded and $x$ is a homogeneous element of $B$.

%For any $n$, if we set $\deg(x) = n$ then $B$ is again connected graded.  
%Clearly $x$ is normal in $B$, and $B/xB \cong A$, so $B$ is a normal extension of $A$.  
%It is well known that $B$ is AS regular if and only if $A$ is.

Suppose that $B$ is any connected $\mb{N}$-graded algebra with a normal non-zero-divisor $x \in B_n$ such that $A \cong B/xB$.  For example, we could have $B = A[x; \sigma]$ where $\sigma$ is a graded automorphism of a connected graded algebra $A$ and $\deg(x) = n$.  The normal extension $B$ of $A$ is again AS regular, with Hilbert series $h_B(t) = h_A(t)/(1-t^n)$ and $\gldim(B) = \gldim(A) + 1$.  Conversely, if $B$ is an AS regular algebra with normal element $x \in B_1$, then $A = B/xB$ will again be AS regular.  In general, however, if $x \in B_n$ is a normal element in an AS regular algebra, then $A =B/xB$ will be AS Gorenstein but usually not regular--this is already clear by taking $x^n$ for $n \geq 2$ in $\kk[x]$.  For proofs see, for example, \cite[Section 1]{RoZh12}.  There are similar statements for the Auslander regular condition \cite{Lev92}.

Le Bruyn, Smith, and Van den Bergh studied central extensions in \cite{LSVdB96}.  Starting with a regular algebra $A$ of dimension 3 with Hilbert series $\frac{1}{(1-t)^3}$, their goal was to classify the possible central extensions $B$ with a central element $x \in B_1$ such that $B/xB \cong A$.  They focused on central extensions rather than the more general normal extensions, because if $x \in B_1$ is normal, there is a Zhang twist $B^{\tau}$ of $B$ such that $x$ becomes central in $B^{\tau}$.  So up to twist equivalence, considering central extensions is sufficient.  For each family of algebras $A$ in the Artin-Schelter classification, the potential relations of central extensions $B$ of $A$ can be written down, and after some reductions, the possible isomorphism types of central extensions were obtained.  There are interesting nontrivial families of central extensions of the $3$-dimensional Sklyanin algebras, for example. 

Cassidy took up this theme in \cite{Cas99}, considering normal extensions of regular algebras of dimension $3$ (including the cubic ones) by normal elements of any degree.  He gave an example of a normal extension by a degree $2$ element which is not twist equivalent to a central extension.  In \cite{Cas03b}, Cassidy pushed the ideas even further to consider normal extensions of those regular algebras $A$ which are not necessarily generated in degree $1$, as had been classified by Stephenson (see Section~\ref{sec:notdeg1} below).

As mentioned above, factoring out a normal element of degree greater than $1$ in a regular algebra does not necessarily yield another regular algebra; it may be only AS Gorenstein.  However, this also means that certain special classes of AS Gorenstein but non-regular algebras $A$ could have normal extensions which are regular.  %In some sense one can see this as the entire strategy of \cite{ATV1}, which classified quadratic regular algebras $A$ of dimension $3$ via their geometric factor rings $A/gA \cong B = B(E, \mc{L}, \sigma)$, where $\deg g = 3$ and $B$ is only AS  Gorenstein.  
In \cite{BaVdB98}, Bauwen and Van den Bergh classify regular algebras $A$ of dimension $4$ with $h_A(t) = \frac{1}{(1-t)^4}$ with a normal element $x \in A_2$ such that $B = A/xA$ is isomorphic to the $2$-Veronese ring $\bigoplus_{n \geq 0} C_{2n}$ of a cubic regular algebra $C$ of dimension $3$; here $B$ is AS Gorenstein only.  Geometrically, these extensions can be seen as noncommutative analogs of the Segre embedding of $\mb{P}^1 \times \mb{P}^1$ into $\mb{P}^3$.

\subsection{Double Ore extensions}

James Zhang and Jun Zhang defined the notion of a double Ore extension in \cite{ZhZh08}.  The general idea is to form a new algebra from a given one by adding two variables both at the same time, with various axioms to ensure that the resulting algebra behaves like a noncommutative version of a polynomial extension in two variables.  Because the relation between the two new variables also needs to be specified, it is a quite complicated object in general.  Here we only define a special case which is the most relevant to the theory of AS regular algebras, following the presentation in \cite{ZhZh09}.

Let $A \subseteq B$ be an extension of connected graded $\kk$-algebras, both generated in degree $1$.  Suppose further that (i) $A_{\geq 1} B = B A_{\geq 1}$; (ii) there is a set $S \subseteq B$ such that $B$ is a free left and right $A$-module with basis $S$; and (iii) the image of $S$ in the quotient ring $C = B/BA_{\geq 1}$ is a $\kk$-basis of $C$. In this setting, we say that $B$ is an \emph{$n$-extension} of $A$ if $C$ is AS regular of dimension $n$, and we say that $B$ is an \emph{$n$-coextension} of $A$ if $A$ is AS regular of dimension $n$.   These notions unify some common constructions already described.  In particular, an extension of connected graded $\kk$-algebras $A \subseteq B$ is a $1$-extension if and only if $B \cong A[x; \sigma, \delta]$ is a graded Ore extension of $A$, and $A \subseteq B$ is a $1$-coextension if and only if $B$ is a normal extension of $A$ \cite[Lemma 1.2]{ZhZh09}.

In this language, we are interested in $2$-extensions $A \subseteq B$.  This is equivalent by \cite[Lemma 1.4]{ZhZh09} to $B$ being a graded double Ore extension of $A$ as defined in \cite{ZhZh08}.  So we assume that $A_{\geq 1} B = B A_{\geq 1}$ and  that $C = B/A_{\geq 1} B$ is regular of dimension $2$.  Without loss of generality, the set $S$ can be taken to have the form $\{ y_1^i y_2^j | i, j \geq 0 \}$ for some elements $y_1, y_2 \in B$, satisfying a relation of the form $y_2 y_1 = p_{12} y_1 y_2 + p_{11} y_1^2 + \tau_1 y_1 + \tau_2 y_2 + y_0$ for some $p_{11}, p_{12} \in \kk$, $\tau_1, \tau_2 \in A_1$, $y_0 \in A_2$, with $p_{12} \neq 0$.  Thus the image $\{ \overline{y_1}^i \overline{y_2}^j | i, j \geq 0 \}$ of this basis in $C$ satisfies the relation $\overline{y_2}\overline{y_1} = p_{12} \overline{y_1}\overline{y_2} + p_{11} \overline{y_1}^2$, which one can check does give $C$ the structure of a regular algebra of dimension $2$.  

One of the advantages of an Ore extension is that given an algebra $A$, we can describe all possible Ore extensions $B = A[x; \sigma, \delta]$ by finding the automorphisms $\sigma$ and $\sigma$-derivations $\delta$ of $A$.  
Describing the $2$-extensions $B$ of $A$ in terms of certain functions of $A$ is part of what is done in \cite{ZhZh08}.  Given the algebra $A$, it requires choosing an algebra homomorphism $\sigma: A \to M_2(A)$ and a $\sigma$-derivation $\delta: A \to A \oplus A$.  Here, we write the elements in the image of $\delta$ as columns, so the usual definition of $\sigma$-derivation requiring $\delta(ab) = \sigma(a) \delta(b) + \delta(a) b$ for all $a, b \in A$ makes sense.  
The multiplication in $B$ is then determined by the matrix equation 
\[
\begin{pmatrix}y_1 \\ y_2  \end{pmatrix} a = \sigma(a)  \begin{pmatrix}y_1 \\ y_2  \end{pmatrix} + \delta(a)
\]
for $a \in A$, together with the relation $y_2y_1 = p_{12} y_1 y_2 + p_{11} y_1^2 + \tau_1 y_1 + \tau_2 y_2 + \tau_0$.
The resulting algebra is written as $A_P[y_1, y_2; \sigma, \delta, \tau]$. 
Unfortunately, in this external description of a double Ore extension there are six additional rather complicated consistency equations involving the data that need to be satisfied, in order to ensure that $A_P[y_1, y_2; \sigma, \delta, \tau]$ is actually a $2$-extension of $A$ \cite[Proposition 1.6]{ZhZh09}.

Double Ore extensions extend the reach of the Ore extension concept enormously, as they turn out to be considerably more general than iterated Ore extensions where one adds $2$ variables, one at a time.  Also, Zhang and Zhang prove that if $B$ is a $2$-extension of $A$, and $A$ is also regular of dimension $n$ (that is, $B$ is also an $n$-coextension of $A$), then $B$ is regular of dimension $n+2$ \cite[Theorem 0.2]{ZhZh08}.  So this really is an effective way to create new examples of regular algebras.  The disadvantage of double Ore extensions is the intricacy of the data defining them, as well as the compatibility conditions required.  Also, it is not known if other basic ring theoretic properties, such as the noetherian condition, are preserved by the construction, unlike the case of usual Ore extensions. 

The most extensive use of the concept so far has been in \cite{ZhZh09}.  In this paper, Zhang and Zhang 
consider double Ore extensions $B = A_P[y_1, y_2; \sigma, \delta, \tau]$ where $A$ is regular of dimension $2$.  Those that are isomorphic to a Ore extension $C[y; \sigma, \delta]$ of a regular algebra $C$ of dimension $3$ are excluded; since regular algebras of dimension $3$ are classified, in theory those Ore extensions can already be understood.  They then completely classify the associated \emph{trimmed} double Ore extensions $A_P[y_1, y_2; \sigma]$ of a regular algebra $A$ of dimension $2$, obtained by setting $\delta = \tau = 0$, obtaining 26 parametrized families. The trimmed double Ore extension $A_P[y_1, y_2; \sigma]$ is an associated graded ring of $A_P[y_1, y_2; \sigma, \delta, \tau]$.  In this particular case, the authors show that all of the trimmed extensions have good ring-theoretic properties (strongly noetherian, Auslander regular, and GK-Cohen Macaulay), and so the same is true by standard arguments for the arbitrary double Ore extensions.  

To give the reader a feel for the concept, here is one family of double Ore extensions from \cite{ZhZh09} (this is family $\mb{C}$ among families labeled $\mb{A}$ through $\mb{Z}$).  Examples like this seem unlikely to have been discovered without the double Ore extension concept, or something similarly involved.
\begin{example} 
\label{ex:double}
The algebra 
\[
\kk \langle x_1, x_2, y_1, y_2 \rangle /(r_1, r_2, r_3, r_4, r_5, r_6)
\]
for $p$ a primitive third root of unity, where the six relations are
\begin{gather*}
r_1 = y_2y_1 - py_1 y_2, \qquad r_2 = x_2x_1 - px_1x_2, \\
r_3 = y_1x_1 + x_1 y_1 - p^2 x_2y_1 - x_1y_2 + p x_2y_2, \\ 
r_4 = y_1x_2  + px_1y_1 - x_2y_1 - x_1y_2 + px_2y_2, \\
r_5 = y_2x_1  + px_1y_1 + 2p^2x_2y_1 - p x_1 y_2 + p x_2 y_2, \\ 
r_6 = y_2x_2  + px_1y_1 - p^2x_2y_1 - x_1y_2 + x_2y_2,
\end{gather*}
is a trimmed double Ore extension $A_P[y_1,y_2; \sigma]$ of $A = \kk \langle x_1, x_2 \rangle/(x_2x_1- p x_1 x_2)$. 
One may observe how the last four relations each mix all of the generators together, unlike what one would expect to see in an iterated Ore extension.  
\end{example}

\subsection{Additional grading}
\label{sec:additional}

While an AS regular algebra is by definition $\mb{N}$-graded, some examples have an additional grading that gives them extra structure.  For example, recall the $S_2$ family of regular algebras of dimension $3$ from Example~\ref{ex:S2}, where
\[
A = A(\alpha) = \kk \langle x_1, x_2, x_3 \rangle/(\alpha x_3x_1 +  x_1 x_3, \ \alpha x_3x_2- x_2x_3, \ x_1^2 - x_2^2).
\]
We see that $A$ is actually $\mb{N} \times \mb{N}$-graded, if we assign $\deg(x_1) = \deg(x_2) = (1,0)$ and $\deg(x_3) = (0,1)$.  The relations are then homogeneous of degrees $(1,1), (1,1)$, and $(2,0)$.
Conversely, of course, any connected $\mb{N} \times \mb{N}$-graded algebra is naturally $\mb{N}$-graded
with $A_n = \bigoplus_{i+j = n} A_{(i,j)}$.  Given a connected $\mb{N}$-graded algebra, generated in degree $1$, we will call it \emph{bigraded} if it has a grading by $\mb{N} \times \mb{N}$ which refines the given one in the sense that $A_n = \bigoplus_{i+j = n} A_{(i,j)}$, and which is nontrivial in the sense that both $A_{(1,0)} \neq 0$ and $A_{(0,1)} \neq 0$.

While the classification of general regular algebras of dimension $4$ remains out of reach, the special class of bigraded regular algebras of dimension $4$ has proven more tractable.  It is easy to see that 
Example~\ref{ex:double} is bigraded by taking $\deg(x_1) = \deg(x_2) = (1,0)$ and $\deg(y_1) = \deg(y_2) = (0,1)$.  In fact, this is a standard feature of a trimmed double extension $A_P[y_1, y_2; \sigma]$.  Conversely, in \cite{ZhZh09}, Zhang and Zhang proved that if $B$ is a regular algebra of dimension $4$ with Hilbert series $\frac{1}{(1-t)^4}$ which is bigraded, then $B$ is either a trimmed double Ore extension $A_P[y_1, y_2; \sigma]$ of a regular algebra $A$ of dimension $2$, or else an Ore extension $A[y; \sigma]$ for a regular algebra $A$ of dimension $3$.  This shows that bigraded regular algebras with this Hilbert series are essentially understood.

In the next sections we will see that bigraded regular algebras of dimension $4$ with other Hilbert series have also been classified.

\section{The Ext-algebra of a regular algebra}
\label{sec:extalgebra}
\subsection{The Ext-algebra}

Let $A$ be any connected $\mb{N}$-graded finitely presented $\kk$-algebra which is generated in degree $1$.  Considering the trivial left module $_A \kk = A/A_{\geq 1}$, we can construct the associated \emph{Ext-algebra} $E = E(A) = \bigoplus_{i=0}^{\infty} \Ext^i_A(\kk,\kk)$.  This has a well-known multiplication given by the Yoneda product, and as such it is again a connected $\mb{N}$-graded $\kk$-algebra.  The study of Ext-algebras and their properties is an important general topic in homological algebra.  In this section we review results which link the AS regularity of a connected graded algebra to properties of its Ext-algebra.

%; for example, one important open question is to find sufficient conditions on $A$ which guarantee that $E(A)$ is finitely generated as an algebra.

\subsection{Koszul algebras}

There is a special case where the connection between a connected graded algebra $A$ and its Ext-algebra $E = E(A)$ is especially tight.  Let 
\[
\dots \to P_n \to \dots \to P_2 \to P_1 \to P_0 \to \kk \to 0
\]
be the minimal graded free resolution of the trivial module $\kk$.  We say that $A$ is \emph{Koszul} if $P_i$ is generated in degree $i$, and so is a direct sum of copies of $A(-i)$, for each $i \geq 0$.  We know that $P_0 = A$, $P_1 = A(-1)^n$ where $n = \dim_{\kk} A_1$ is the minimal number of generators of $A$, and $P_2 = \bigoplus_{i=1}^s A(-d_i)$, where the minimal homogeneous relations are $r_1, \dots, r_s$ with $\deg r_i = d_i$.  In particular, if $A$ is Koszul then it is quadratic, that is, its minimal relations all have degree $2$.  

For any connected $\mb{N}$-graded $\kk$-algebra $A$ generated in degree $1$ with quadratic relations, say $A = \kk \langle x_1, \dots, x_n \rangle/(r_1, \dots, r_m)$, with $\deg r_i = 2$ for all $i$, the \emph{Koszul dual} of $A$ is an algebra $A^!$ defined as follows.  Let $V = A_1$ with its $\kk$-basis $x_1, \dots, x_n$, and let $x_1^*, \dots, x_n^*$ be a dual basis of $V^*$.  Thinking of $r_1, \dots, r_m$ as elements of $V \otimes_{\kk} V$, say with $\kk$-span $R$, we let $R^{\perp} \in V^* \otimes V^*$ be the orthogonal complement of $R$ under the canonical pairing between $V^* \otimes V^*$ and $V \otimes V$.  Then we define $A^! = \kk \langle x_1^*, \dots, x_n^* \rangle/(R^{\perp})$.  For example, it is easy to calculate that if $A = \kk [x_1, \dots, x_n]$ is a commutative polynomial ring, then $A^!$ is isomorphic to an exterior algebra in $n$ generators.  If $A$ is Koszul then so is $A^!$.  A basic result is that a connected graded $\kk$-algebra $A$, generated in degree $1$ with quadratic relations, is Koszul if and only if $A^! \cong E(A) = \bigoplus_{i = 0}^{\infty} \Ext^i_A(k,k)$ as graded $\kk$-algebras.    Note that by definition $(A^!)^! \cong A$ always, so one consequence of this is that if $A$ is Koszul, then $E(E(A)) \cong A$.

Given an AS regular algebra, whether or not it is Koszul can generally be easily recognized by its Hilbert series.  If $A$ is AS regular of global dimension $d$ with Hilbert series $h_A(t) = \frac{1}{(1-t)^d}$ (the same as a polynomial ring in $d$ variables of degree $1$), then in fact $A$ is Koszul, as noted by Zhang \cite[Theorem 5.11]{Smi96}.  Thus among regular algebras of dimension $3$, for example, the quadratic ones are Koszul, and the cubic ones are not, because Koszul algebras must have quadratic relations.  It seems very likely that, conversely, any Koszul AS regular algebra of dimension $d$ must have Hilbert series $\frac{1}{(1-t)^d}$, but this seems to be unknown in general.  This is true when $d \leq 5$ under mild additional hypotheses, by the classification of possible Hilbert series we will review in Section~\ref{sec:hilb} below.

Suppose now that $A$ is an AS regular algebra, not necessarily Koszul.  Since $\gldim(A) = d < \infty$, we have $E = E(A) = \bigoplus_{i=0}^d \Ext^i_A(\kk,\kk)$.  Even if we do not necessarily assume that $A$ is noetherian, for a regular algebra each $P_i$ in the minimal graded free resolution of $\kk$ is finitely generated \cite[Propostion 3.1]{StZh97}, and this implies that 
$E_i = \Ext^i_A(\kk, \kk)$ is finite-dimensional over $\kk$ for all $i$.    So $E$ is a finite-dimensional connected $\mb{N}$-graded algebra.  A finite dimensional algebra $E$ is called \emph{Frobenius} if $E^* = \Hom_{\kk}(E, \kk)$ is isomorphic to $E$ as a left $E$-module.

We can now see how the regularity of an algebra is reflected in its Ext-algebra.
\begin{theorem}\cite[Corollary D]{LPWZ08} 
\label{thm:frobenius}
Let $A$ be a connected $\mb{N}$-graded algebra, generated in degree $1$, with $\gldim(A) < \infty$. 
Then $A$ satisfies the AS Gorenstein condition if and only if $E(A)$ is Frobenius.  In particular, if $\GKdim(A) < \infty$ then $A$ is regular if and only if $E(A)$ is Frobenius.
\end{theorem}
We have stated the final, most general version of this theorem, but important special cases of the theorem were proved earlier.  In particular, Smith considered the case where $A$ is Koszul and proved Theorem~\ref{thm:frobenius} in this case \cite{Smi96}.  Inspired partly by a desire to put cubic regular algebras of dimension $3$ on a similar footing as the quadratic ones, Berger defined a notion of \emph{$m$-Koszul} algebra \cite{Ber01}.  This is an algebra $A$ presented by homogeneous relations, all of some degree $m$, for which the term $P_i$ in the minimal graded free resolution of $\kk$ is generated in degree $d_i$ for all $i$, where for all even $i$ we have $d_i = mi/2$ and $d_{i+1} = (mi/2) +1$.  Thus Koszul algebras are $2$-Koszul, and cubic regular algebras of dimension $3$ are $3$-Koszul.  Berger and Marconnet  showed that Theorem~\ref{thm:frobenius} holds for $m$-Koszul algebras \cite[Theorem 1.2]{BerMa06}.  In addition to cubic regular algebras this applies to some examples with infinite GK-dimension (see also Section~\ref{sec:weak} below).

We are unaware of any explicit examples of AS regular algebras which are $m$-Koszul for $m > 2$, other than the cubic regular algebras of dimension $3$; in fact, non-quadratic regular algebras of higher dimension seldom have relations all in one degree, as we will see in Section~\ref{sec:Hilbert}.  Another generalization of Koszul, which subsumes $m$-Koszulity, is the notion of $\mc{K}_2$-algebra defined by Cassidy and Shelton, which does hold for certain regular algebras of dimension $4$ with relations in multiple degrees \cite[Theorem 4.7]{CaSh08}.

In order to obtain the full generality of Theorem~\ref{thm:frobenius}, Lu, Palmieri, Wu, and Zhang introduced new methods, which we describe next.

\subsection{$A_{\infty}$-algebras}
\label{sec:Ainfinity}

$A_{\infty}$-algebras were defined by Stasheff.  The original motivation for their introduction came from topology, but as often happens, they have found many applications in an algebraic setting far removed from topology.
\begin{definition}
    An $A_{\infty}$-algebra is a $\mb{Z}$-graded vector space $\bigoplus_{i \in \mb{Z}} A^i$ together with a family of $\kk$-linear ``multiplication maps" $m_n: A^{\otimes n} \to A$ for $n \geq 1$, where $m_n$ is graded of degree $2-n$, such that Stasheff's identities are satisfied: 
    \[
\sum (-1)^{r+st} m_{r + 1 + t}(\on{id}^{r} \otimes m_s \otimes \on{id}^t) = 0
    \]
for all $n \geq 1$, where we sum over all $n = r + s + t$ with $r, t \geq 0$ and $s \geq 1$.    
\end{definition}
We give the definition here explicitly in order to demonstrate that the relations that must hold among the $m_n$ are quite subtle.  The map $m_2: A \otimes A \to A$ is a degree $0$ multiplication operation, while $m_1: A \to A$ is of degree $-1$ and is a differential of $\bigoplus_{i \in \mb{Z}} A^i$ considered as a complex.  In fact, when $m_n = 0$ for $n \geq 3$, an $A_{\infty}$-algebra is precisely a differential graded algebra.  The maps $m_n$ for $n \geq 3$ are thought of as ``higher multiplications", and in general the multiplication $m_2$ is associative only up to homotopy (although it is actually associative in case $m_1 = 0$ or $m_3 = 0$).  For a general introduction to $A_{\infty}$-algebras see \cite{Ke01}, and for an introduction geared towards their application to noncommutative algebra see \cite{LPWZ04}.  

Suppose now that $A$ is a connected $\mb{N}$-graded $\kk$-algebra, finitely generated in degree $1$.
Then the Ext-algebra $E = E(A) = \bigoplus_{i \geq 0} \Ext^i_A(\kk, \kk)$ can be given the structure of an $A_{\infty}$-algebra for which $m_1 = 0$ (such $A_{\infty}$-algebras with $0$ differential $m_1$ are called \emph{minimal}).  Moreover, the relations of the algebra $A$ can be recovered from knowledge of $E$ together with its entire $A_{\infty}$-algebra structure, in a straightforward way \cite[Corollary B]{LPWZ09}.  If $A$ is Koszul, then in fact $m_n = 0$ for $n \neq 2$, so $E$ is just an associative algebra and this is just the statement that $A$ is recoverable from the algebra $E$ (via $E(E(A)) \cong A$).  If $A$ is not Koszul, on the other hand, one needs to remember the higher multiplications in order to recover $A$ from $E$.

Using these ideas, Lu, Palmieri,  Wu, and Zhang proved Theorem~\ref{thm:frobenius} in its full stated generality \cite[Corollary D]{LPWZ08}.  Thus, for any connected graded algebra of finite global dimension, the Gorenstein condition is equivalent to the condition that the Ext-algebra $E(A)$ is Frobenius.  Because one cannot recover $A$ from $E = E(A)$ without invoking the $A_{\infty}$-algebra structure on $E$ in general, it is not surprising that one needs $A_{\infty}$-techniques to characterize when $A$ is Gorenstein in terms of the properties of $E$, even though in the end the Frobenius property of $E$ does not depend on the higher multiplications.

\begin{example}[Example 13.4, \cite{LPWZ04}]
Consider the algebras 
\[
A = A(\alpha, \beta) = \kk \langle x, y \rangle /(xy^2 + \alpha y^2x + \beta yxy, \ x^2y + \alpha yx^2 + \beta xyx).
\]
with $\alpha, \beta \in \kk$.  These are the graded algebras among the well-studied class of \emph{down-up algebras} which 
were first defined by Benkart and Roby.

The Ext-algebra $E = E(A)$ has Hilbert series $h_E(t) = 1 + 2t + 2t^2 + t^3$; choosing a basis as follows:
\[
E = E^0 \oplus E^1 \oplus E^2 \oplus E^3 = \kk1 \oplus (\kk a \oplus \kk b) \oplus (\kk c \oplus \kk d ) \oplus \kk e, 
\]
where $1$ is the identity, then the multiplication is given by the rules $ac = e$, $bd = \alpha e$, $ca = \alpha e$, $db = e$, with all other products of pairs of basis elements equal to $0$.  It is easy to see directly that $E$ is Frobenius exactly when $\alpha \neq 0$, so $A$ is regular of dimension $3$ under the same condition.

The higher multiplications on $E$, when restricted to elements in $E^1$, are determined from the relations in a straightforward way \cite[Corollary B]{LPWZ09}.  In this case $m_n = 0$ for $n \geq 4$, 
and $m_3(a, b, b) = c$, $m_3(b, b, a) = \alpha c$, $m_3(b, a, b) = \beta c$, $m_3(a,a,b) = d$, 
$m_3(b, a, a) = \alpha d$, $m_3(a, b, a) = \beta d$, and $m_3(a,a,a) = m_3(b,b,b) = 0$.

Note that the multiplication $m_2$ of the Ext-algebra $E$ only involves $\alpha$, while $\beta$ does not appear.  As $\beta$ varies, one gets a whole family of regular algebras with isomorphic Ext-algebras.  On the other hand, $\beta$ appears in the formulas for the higher multiplication $m_3$, which is why $A$ can be recovered from $E$ with its entire $A_{\infty}$-structure.  
\end{example}

\subsection{Using the Ext-algebra}

One of the main motivations in relating properties of a graded ring and its Ext-algebra, as in Theorem~\ref{thm:frobenius}, is the hope that one might be able to understand or classify regular algebras by instead considering finite-dimensional Frobenius algebras.  There has been some success in this regard.

In \cite{LPWZ07}, the authors apply this idea to regular algebras $A$ of dimension $4$, in particular those with Hilbert series $h_A(t) = \frac{1}{(1-t)^2(1-t^2)(1-t^3)}$.  These algebras have $2$ degree $1$ generators and $2$ minimal relations, one each of degrees 3 and 4. We discuss the possible Hilbert series of regular algebras of dimension $4$ in the next section; the case at hand is the one in which the Ext-algebra $E = E(A)$ is smallest, specifically $h_E(t) = 1 + 2t + 2t^2 + 2t^3 + t^4$.
%, and they refer to this case as type $(12221)$.  
The authors show that the only nonzero higher multiplications of $E$ are $m_2, m_3, m_4$.  Starting with a Frobenius algebra $E$ of this shape with unknown parameters and arbitrary $m_3, m_4$ satisfying the Shasheff relations, they classify those that are ``generic" in the sense that the maps $m_3, m_4$  satisfy certain Zariski open conditions.  Then, relations for the corresponding algebras $A$ such that $E = E(A)$ can be determined.  Finally, which of these $A$ are actually noetherian and regular is decided.  The analysis is quite sensitive, but leads in the end to just $4$ parametrized families of regular algebras \cite[Theorem A]{LPWZ07}.  This classifies all noetherian regular algebras of dimension $4$ with this Hilbert series whose Ext-algebra has generic higher multiplications.  

As a sample of what these algebras can look like, here is one of the families: 
\[
A(p) = \kk \langle x, y \rangle/(xy^2 - p^2 y^2x, \ x^3y + p x^2yx + p^2 xyx^2 + p^3 yx^3),\ \qquad \text{where}\ 0 \neq p \in \kk.
\]
For this family, $A(p)$ is indeed regular for all nonzero $p$ (in fact, all algebras in the family are Zhang twists of $A(1)$).   

Note that all of the members of this family $A(p)$ are bigraded with $\deg(x) = (1,0)$ and $\deg(y) = (0,1)$.  This is not an accident, as it turns out that all of the algebras classified in \cite[Theorem A]{LPWZ07} are bigraded.  The authors prove that the same four families also classify regular algebras of dimension $4$ with Hilbert series $\frac{1}{(1-t)^2(1-t^2)(1-t^3)}$ that are bigraded \cite[Theorem B]{LPWZ07}.
  
More recent work has pushed this $A_{\infty}$-algebra method a bit farther.  Many authors have studied PBW deformations of regular algebras---these are algebras which have a nice filtration such that the associated graded algebra is AS regular.  Shen, Zhou, and Lu exploit this idea to construct new regular algebras closely related to given ones, produced by ``adding terms" to the relations \cite{SZL15} \cite{SZL16}.  This allows them to extend the technique of $A_{\infty}$-algebras from \cite{LPWZ07}  to classify certain regular algebras with the same Hilbert series where the $A_{\infty}$-structure on the Ext-algebra is non-generic.

In another direct generalization of the work of \cite{LPWZ07}, Wang and Wu  looked at regular algebras of dimension $5$ generated by two degree $1$ generators, with relations of degrees $(4, 4, 4)$ \cite{WangWu12}, which have Hilbert series $\frac{1}{(1-t)^2(1-t^2)(1-t^3)^2}$.  They performed an analysis similar to \cite{LPWZ07}, and classified all such algebras which are domains, under a genericity condition on the structure of the associated $A_{\infty}$-algebra.  Similarly as in \cite{LPWZ07}, the algebras that arise under this generic condition are also bigraded.
%For a few more papers that exploit the Koszul dual in various ways to study the properties of a regular algebra or construct new regular algebras, see \cite{ShTi01}, \cite{Mor06a}.
%For a Koszul regular algebra $A$ of dimension $4$, with Hilbert series $1/(1-t)^4$, its Frobenius Ext algebra $E(A)$ has Hilbert series $1 + 4t + 6t^2 + 4t^3 + t^4$.  This is quite a large algebra, and classifying Frobenius algebras of this shape appears very difficult.  There are a few papers we can mention that exploit different kinds of relations between $A$ and its Koszul dual.  Shelton and Tingey 
%\cite{ShTi01} showed that by starting with such a Koszul regular algebra $A$, if $A$ has a sequence of elements $y_1, y_2, y_3, y_4 \in A_2$ which is regular and normal (that is $y_i$ is a normal non-zero-divisor in the ring $A/(y_1, \dots, y_{i-1})$ for each $i$), then $E = A/(y_1, \dots, y_4)$ has Hilbert series $1 + 4t + 6t^2 + 4t^3 + t^4$ and $E^!$ is another regular algebra of dimension $4$.  

\section{Hilbert series}
\label{sec:Hilbert}
Let $A$ be a regular algebra, generated in degree $1$, with global dimension $d$.  Much of the work on regular algebras we described above has concentrated on those which are Koszul with $h_A(t) = \frac{1}{(1-t)^d}$, and many of the most prominent examples, such as the Sklyanin algebras, are in this category.  But even for $d = 3$, in addition to the Koszul quadratic regular algebras there are also the cubic regular algebras with Hilbert series $h_A(t) = \frac{1}{(1-t)^2(1-t^2)}$, and the work of Artin, Tate, and Van den Bergh managed to treat both cases in roughly the same way.

For higher $d$, an obvious general question is what Hilbert series are possible for regular algebras.  We will see momentarily that there are not many possibilities for $d = 4$, but starting with $d = 5$ there are a considerable number of cases.  It seems unlikely that any eventual understanding or classification will be able to treat these all uniformly.

Already in \cite{AS87}, Artin and Schelter gave examples of regular algebras of dimension $4$ 
%which are universal enveloping algebras of graded Lie algebras, 
with Hilbert series $\frac{1}{(1-t)^4}$, $\frac{1}{(1-t)^3(1-t^2)}$, and $\frac{1}{(1-t)^2(1-t^2)(1-t^3)}$.  Lu, Palmieri, Wu, and Zhang proved in \cite{LPWZ07} that if $A$ is regular, generated in degree $1$, with $\on{gl.dim}(A) = 4$, and either $\GK(A) \geq 3$ or $A$ is a domain, then $A$ has one of these three Hilbert series (and in fact $\GK(A) = 4$).  Of course, conjecturally a regular algebra $A$ is a domain and satisfies $\GK(A) = \gldim(A)$ (as in Question~\ref{ques:props}).  So the hypothesis that $\GK(A) \geq 3$ or $A$ is a domain is almost surely unnecessary, and one expects that these three Hilbert series are the only possible ones in dimension $4$.

\subsection{Enveloping Algebras and graded Ore extensions}

The easiest way to construct examples of higher-dimensional regular algebras with various Hilbert series is to use graded Ore extensions.  If $A$ is AS regular, then an Ore extension $B = A[x; \sigma, \delta]$ is again connected graded and regular by taking $\deg(x) = n$, as long as $\sigma$ preserves degree and $\delta$ is a $\sigma$-derivation of degree $n$.
We are still primarily interested in regular algebras which are generated in degree $1$; the key observation is that $B$ may be generated in degree $1$ even if $A$ is not.

For instance, let $L$ be a finite-dimensional Lie algebra over $\kk$.  Assume that $L$ is positively graded, so $L = \bigoplus_{i=1}^m L_i$ 
as vector spaces, with $[L_i, L_j] \subseteq L_{i+j}$ for all $i,j$.  Let $A = U(L)$ be the universal enveloping algebra of $L$.  Then because $L$ is graded, $A$ is naturally connected $\mb{N}$-graded.  Under these conditions $U(L)$ is Artin-Schelter regular, with Hilbert series $h_{U(L)}(t) =  \frac{1}{\prod_{i=1}^m (1-t^i)^{\dim_{\kk} L_i}}$.  If $L$ is generated as a Lie algebra by $L_1$, then $U(L)$ will be generated as an algebra in degree $1$.  The algebra $U(L)$ can always be expressed as an iterated Ore extension of the form $U(L) = \kk[x_1][x_2; \delta_2] \dots [x_m; \delta_m]$, by choosing $x_i$ to be a homogeneous basis of $L$ in decreasing order of degree.  

\begin{example} 
Let $L$ have $\kk$-basis $x, y, w, z$ with $\deg(x) = \deg(y) = 1$, $\deg(w) = 2$, and $\deg(z) = 3$, where $[x, y] = w = -[y,x]$, $[x, w] = z = -[w,x]$, and all other brackets between pairs of variables $0$.  This is easily seen to be a graded Lie algebra generated in degree $1$.  Then $U(L)$ is regular with Hilbert series $\frac{1}{(1-t)^2(1-t^2)(1-t^3)}$, generated by only $x$ and $y$ because of the relations $w = xy-yx$ and $z = xw-wx$.  We can also write $U(L) = \kk[z][w][y][x; \delta]$ where 
$\delta(y) = w$, $\delta(w) = z$, and $\delta(z) = 0$.
\end{example}

As an aside, we mention that there is also a different way to use Lie algebras to construct interesting examples of regular algebras.
When $L$ is not a graded Lie algebra, $U(L)$ will not be graded either.  However, in this case one can form a \emph{homogenized enveloping algebra} by adding an additional variable $z$ and defining 
\[H(L) = \kk \langle x_1, \dots, x_n, z \rangle/(x_jx_i - x_ix_j - [x_i, x_j]z \,| \, 1 \leq i < j \leq n)
\]
for a basis $x_1, \dots, x_n$ of $L$.  These algebras need not be iterated Ore extensions.  Le Bruyn and Van den Bergh proved that homogenized enveloping algebras $H(L)$ are AS regular, and they gave some information about their linear modules \cite{LeBVdB93}.  Le Bruyn and Smith further investigated $H(L)$ for the special case of $L = \mf{sl}_2$ in \cite{LeBSm93}.

\subsection{Hilbert series in dimension 5}

For each of the three possible Hilbert series for regular algebras of dimension $4$, there is a graded Lie algebra $L$ of dimension $4$ such that $U(L)$ is regular (generated in degree $1$) of that Hilbert series.  In fact the examples initially given by Artin and Schelter are of that form \cite[Example 1.21]{AS87}.

This raises the question of whether universal enveloping algebras already achieve all of the possible Hilbert series of regular algebras in higher dimension.   Actually, already in global dimension $5$ there is an additional consideration; one may ask for the finer information of the possible \emph{Betti numbers} of the minimal graded free resolution of $\kk$.  That is, one wants to know the non-negative integers $\beta_{ij}$ such that in the minimal graded resolution $0 \to P_5 \to \dots \to P_2 \to P_1 \to P_0 \to \kk$, one has $P_i = \bigoplus_{j=1}^{d_j} A(-\beta_{ij})$.  The Hilbert series can be calculated from the Betti numbers, but starting in global dimension $5$ it is possible for regular algebras to have different Betti numbers but the same Hilbert series.  Assuming we have $m$ degree $1$ generators and minimal relations $r_1, \dots, r_t$ of degrees $\deg(r_i) = d_i$, then $P_1 = A(-1)^m$ and $P_2 = \bigoplus_{i=1}^r A(-d_i)$.  By the symmetry of the free resolution coming from the Gorenstein condition, there is an $\ell$ such that $P_5 = A(-\ell)$ (the number $\ell$ is sometimes called the \emph{AS-index} or the \emph{Gorenstein parameter}). Once $\ell$ is fixed 
the free resolution must look like
\[
0 \to A(-\ell) \to A(-\ell+1)^m \to \bigoplus_{i=1}^r A(-\ell+d_i) 
\to \bigoplus_{i=1}^r A(-d_i) \to A(-1)^m \to A \to \kk \to 0,
\]
so the Betti numbers are determined.  Then the Hilbert series is equal to $h_A(t) = \frac{1}{p(t)}$ where 
$p(t) = 1 - mt + \sum_{i=1}^r t^{d_i} - \sum_{i=1}^r t^{\ell-d_i} + mt^{\ell-1} - t^{\ell}$.

Floystad and Vatne considered the question of the possible Hilbert series of dimension $5$ regular algebras $A$ (generated in degree $1$).  It is relatively straightforward to calculate all of the Hilbert series of universal enveloping algebras $U(L)$ where $L$ is a graded Lie algebra of dimension $5$ generated in degree $1$, which is done in \cite[Proposition 3.1]{FlVa11}.  The authors then 
proved that if $A$ is regular of global dimension $5$, generated by $m = 2$ degree $1$ elements, has $\GK(A) \geq 4$, and $A$ is a domain, then the only possibilities for the degrees $(d_1, \dots, d_t)$ of the minimal relations are $(4,4, 4)$, $(4,4,4, 5)$,  $(4, 4,4,5, 5)$, $(3, 5, 5)$, or $(3, 4, 7)$.  In each case there is only one possible $\ell$ (respectively $\ell = 10, 10, 10, 11$ or $12$) \cite[Theorem 5.6]{FlVa11}.  The respective Hilbert series are $\frac{1}{(1-t)^2(1-t^2)(1-t^3)(1-t^3)}$ (for the first three possibilities), $\frac{1}{(1-t)^2(1-t^2)(1-t^3)(1-t^4)}$, and $\frac{1}{(1-t)^2(1-t^2)(1-t^3)(1-t^5)}$.   Note that certainly all of the Hilbert series which occur here are the same as a commutative graded polynomial ring with variables of some degrees, consistent with a positive answer to Question~\ref{ques:props}(5). 

The authors gave examples of regular algebras with each of the types with three relations (leaving open the question of whether any regular algebras with two generators and minimal relations of degrees $(4,4,4, 5)$ or $(4, 4, 4, 5, 5)$ exist). They refer to the case with relations of degree $(3, 4, 7)$  as ``extremal" since it has the largest AS-index of $12$.  From their earlier calculation of the possible Hilbert series of universal enveloping algebras $U(L)$, they conclude that there is no enveloping algebra of this extremal type.  Thus enveloping algebras of graded Lie algebras are not in fact general enough to achieve all possible Hilbert series of regular algebras of dimension $5$.

Since enveloping algebras of positively graded Lie algebras are just a special case of iterated Ore extensions, one may also ask whether graded iterated Ore extensions exhibit all possible types of regular algebras of dimension $5$.  Elle studied iterated graded Ore extensions of dimension $5$ which are generated in degree $1$ \cite{Elle17}.  She showed that for all types in dimension $5$ for which there is some known example there is an example of that type given by an iterated Ore extension.  In particular, this is true for the extremal type which Floystad and Vatne showed cannot be obtained via a universal enveloping algebra of a graded Lie algebra.  She also exhibited a similar type with three generators and relations of degrees $(2, 2, 3)$ which can be obtained via an iterated Ore extension but not via an enveloping algebra.  There is still no known example to our knowledge of a Hilbert series achieved by some regular algebra but not by a graded iterated Ore extension.

This type of analysis becomes harder as the dimension $d$ of the regular algebra increases.  Relatively little is known about the full range of possible Hilbert series and Betti numbers for regular algebras of dimension $d$ with $d \geq 6$. It would be good to have a more theoretical understanding of the possibilities, apart from brute force computation.

\subsection{Gr\"obner basis methods}

The original work of Artin and Schelter made heavy use of the methods of noncommutative Gr\"obner bases to construct examples of regular algebras of dimension $3$ \cite[Sections 5-10]{AS87}. 
For certain special kinds of regular algebras of dimension $4$ and $5$, Gr\"obner basis methods are in fact sufficient to achieve a classification.  

Recall that bigraded regular algebras of dimension $4$ with Hilbert series $\frac{1}{(1-t)^4}$ were classified in \cite{ZhZh09} by Zhang and Zhang using double Ore extensions, while bigraded regular algebras with Hilbert series $\frac{1}{(1-t)^2(1-t^2)(1-t^3)}$ were classified by Lu, Palmieri, Wu, and Zhang by using $A_{\infty}$-algebras \cite{LPWZ07}.   That leaves algebras with Hilbert series $\frac{1}{(1-t)^3(1-t^2)}$, which have $3$ degree $1$ generators, and $4$ relations, $2$ each of degrees $2$ and $3$.  The author and Zhang considered bigraded regular algebras of this type in \cite{RoZh12}. It turns out that after some reductions, including eliminating Ore extensions of regular algebras of dimension $3$, all such algebras can be presented with relations with the same leading terms, and one can show that for the Hilbert series to be correct, the Gr\"obner basis must have a specific form (namely, the overlaps must lead to just one additional relation with a particular leading term, after which all overlaps must resolve).  One can then use a computer calculation to find families of algebras for which the Gr\"obner basis has this form.  In \cite{RoZh12} the authors show that such examples fall into one of $8$ parameterized families.  Computational methods are also used to help check that the algebras in these families really are regular.  

Thus  bigraded algebras of dimension $4$ of all three possible Hilbert series are understood.  One goal of classifying algebras with additional special properties is to produce a large portfolio of examples on which to test the potential validity of conjectures about regular algebras.  Using the various classification results, all bigraded regular algebras of dimension $4$ are known to be strongly noetherian, Auslander regular, and GK-Cohen Macaulay, for example, giving further evidence that these properties might hold for all regular algebras. 

Zhou and Lu \cite{ZhouLu14}  classified all nontrivially bigraded regular algebras of dimension $5$ with two generators, using a Gr\"obner basis method similar to that of \cite{RoZh12}.  There are such examples of most of the possible relation types $(4, 4, 4)$, $(3, 5, 5)$, $(3, 4, 7)$ and $(4, 4, 4, 5, 5)$ in Floystad and Vatne's work.  In particular, these were the first examples of the final relation type, but it is still open if the relation type $(4, 4, 4, 5)$ can occur.  The authors also note that in all of their examples, there is a choice of Gr\"obner basis for the ideal of relations whose leading terms consist of \emph{Lyndon words}: words $w$ that are maximal in the lexicographical ordering among all of the cyclic permutations of their letters.  A Gr\"obner basis consisting of Lyndon words has good combinatorial properties, and so the authors ask if for bigraded regular algebras of arbitrary dimension with 2 generators, it is always true that there is a Gr\"obner basis of Lyndon words.  This would imply several of the conjectured properties of regular algebras in this special case, for example it would follow that $\GK(A) = \gldim(A)$, and that the Hilbert series of $A$ is the same as that of some graded polynomial ring \cite[Theorem 1.9]{ZhouLu14}.  

In higher dimension, Gateva-Ivanova wrote down all regular algebras occurring as universal enveloping algebras of certain kinds of Lie algebras of dimension $6$ and $7$---those with a Gr\"obner basis of relations with leading terms consisting of Lyndon words \cite{G-I22}.

\section{Generalizations of the AS regular condition}
In this section we discuss what happens when some of the hypotheses in the definition of AS regular are relaxed.  Then in the final section we will return to some of the recent work on AS regular algebras.

\subsection{Generators in arbitrary degrees}
\label{sec:notdeg1}
Artin and Schelter did not require generation in degree $1$ in their original definition of regular algebra, but almost all of the first papers in the subject only considered regular algebras under this additional assumption, as we have in this survey until now.  It was clear from the outset that allowing generators in arbitrary degrees would introduce many significant complications.  Recall that point modules, their parametrizing scheme, and the associated automorphism coming from the truncation-shift operation played a large role in the classification of regular algebras of dimension $3$.  If $A$ is a connected $\mb{N}$-graded $\kk$-algebra with generators in arbitrary degrees, and $M$ is a point module for $A$ as defined above, then there is no reason the truncation shift $M_{\geq 1}(1)$ must be a point module any longer---it still has the correct Hilbert series, but it may not be cyclic.  In \cite{StaZh94a}, Stafford and Zhang presented an interesting example of a (non-regular) idealizer ring $B$, which has generators of degree $1$ and $2$ inside of a Jordan plane.  This ring has many exotic properties in the context of noncommutative projective geometry. In particular, it has point modules $M$ for which $M_{\geq 1}(1)$ is not a point module.

The basic analysis of how the axioms of a regular algebra restrict the form of the free resolution of $\kk$, as we outlined earlier in Section~\ref{sec:AS}, does go through with minor adjustments for algebras with generators in arbitrary degrees.  For instance, the classification of regular algebras of dimension $2$ generated in arbitrary degrees is very similar:
\begin{proposition}[Proposition 3.3, Remark 3.4, \cite{StZh97}]
Let $A$ be an AS regular algebra of global dimension $2$.  Then either 
\begin{enumerate}
    \item $A \cong \kk \langle x, y \rangle/(yx-qxy)$ for some $0 \neq q \in \kk$, where $\deg x = a$, $\deg y = b$, for any 
    $a, b \geq 1$; or
    \item $A \cong \kk \langle x, y \rangle/(yx - xy - x^m)$ where $\deg x = a$, $\deg y = (m-1)a$, for some $a \geq 1$ and $m \geq 2$.
\end{enumerate}
\end{proposition}

In fact, many aspects of the classification of regular algebras of dimension $3$ also go through for generators in arbitrary degrees.  The theory was worked out by Stephenson in a series of papers \cite{Ste96} \cite{Ste97}.  In particular, it is still possible to use point module techniques, but this requires care because of the issues mentioned above.  For these regular algebras, it turns out that if one considers only point modules which are generated in degrees $0$ and $1$ and have no socle, then there is a nice scheme $X$ parametrizing them for which the methods of \cite{ATV1} can then be used. For instance, there are Sklyanin-like examples for which $X$ is a smooth elliptic curve.  %The methods are particular to these algebras, and it remains unclear if there is a ``best" way to define point modules which works uniformly for general algebras generated in arbitrary degrees.

Stephenson also investigated the interesting geometry associated with these regular algebras.
The particular regular algebras $A$ with $3$ generators in degrees $(1, 1, n)$ for some $n \geq 2$ are examined more closely in \cite{Ste00} \cite{Ste00corr}.  The associated category $\rproj A$ can be thought of as a noncommutative analog of a weighted projective space $Z = \on{Proj} \, k[x, y, z]$, where $\deg(x) = \deg(y) = 1$, and $\deg(z) = n$.  The commutative geometry of $Z$ is well-known: $Z$ can be obtained by blowing down the unique exceptional curve on a rational ruled surface with invariant $n$.  In terms of the point module techniques used in \cite{Ste97}, Stephenson shows that the point modules for $A$ are parametrized by $X$ where either $X \cong Z$ or else $X$ is isomorphic to a subscheme of $Z$ defined by a single equation of degree $n +2$.  In the latter case, $A$ has a normal element $g$ of degree $n + 2$ such that $A/gA$ is isomorphic to a twisted homogeneous coordinate ring on $X$.  Thus, for these algebras the analogy with the geometric approach used in \cite{ATV1} is very tight.  Ingalls and Patrick \cite{InPa02} made the geometric analogy even more explicit in terms of the noncommutative projective schemes, by showing there is a  way to define a blowup of $\rProj A$ at a point using a Rees ring, such that the resulting category is a quantum ruled surface in a sense studied earlier by Patrick \cite{Pat97} \cite{Pat00}.  

Smith studied another regular algebra $A$ of dimension $3$ with two generators of degrees $1$ and $2$, and showed that $\rproj A$ is equivalent to the category of coherent sheaves on a degree $6$ Del Pezzo surface \cite{Smi12b}.

One important reason why one might be interested in regular algebras generated in arbitrary degrees is their occurrence as invariant rings in noncommutative invariant theory (see \cite{Kir16}).
%We are unaware of any systematic study of such regular algebras of higher dimension than $3$, however.

%If $G$ is a finite group of graded automorphisms of a regular algebra $A$, then one is interested in understanding the properties of the invariant subring $A^G = \{ a \in A | g(a) = a\ \text{for all}\ g \in G \}$, for example under what circumstances it might again be AS regular or AS Gorenstein.  One might well take $A$ to be generated in degree $1$, but of course the ring $A^G$ is almost never generated in degree $1$.  Indeed in commutative invariant theory, where $A$ is a polynomial ring, the degrees in which the minimal generators of $A^G$ lie is one of the important questions in the theory.  Thus examples of regular algebras with generators in arbitrary degrees appear incidentally as invariant rings in the literature of noncommutative invariant theory.  However, we are unaware of any systematic study of regular algebras with arbitrary generating degrees since Stephenson's work.

\subsection{Omitting the condition on growth}
\label{sec:weak}

Consider the definition of an AS regular algebra $A$, but remove axiom (2) that the graded pieces of $A$ have polynomially bounded growth (or equivalently, that $\GK(A)$ is finite).  We call such an algebra \emph{weakly} AS regular.  There is no standard term; some papers in the literature simply call such examples AS regular again, making it clear at the outset that they are considering a version of the definition in which finite GK dimension is not required, unlike Artin and Schelter's original definition.

Suppose that $A$ is weakly regular (and not regular), so that $\GK(A) = \infty$.  As we remarked earlier in Section~\ref{sec:GK}, in this case $A$ has exponential growth, and therefore cannot be noetherian by a theorem of Stephenson and Zhang \cite{StZh97}.  Thus a weakly regular algebra certainly does not closely resemble a commutative polynomial ring, and is a different kind of animal.  Techniques coming from point modules are also likely to be difficult to apply: one does not expect the inverse system of truncated point modules defining the point scheme to stabilize in a non-noetherian algebra.  Nevertheless, weakly regular algebras are interesting and seem to be well-behaved in certain other ways.  For example, it is still true that all known examples are domains.

Once again, the examples of global dimension $2$ are special enough to be classified, as was done by Zhang.  For simplicity 
we state the result for those algebras generated in degree $1$, though the general classification is not really any harder. 
\begin{theorem} \cite{Zha98a}
\label{thm:zhang}
Let $A$ be a weakly AS regular algebra of global dimension $2$ which is generated in degree $1$.  Then 
\[
A \cong \kk \langle x_1, \dots, x_n \rangle/(\sum_{i=1}^n x_i \tau(x_i) )
\]
For some $\kk$-linear bijection $\tau$ from $\kk x_1 + \dots + \kk x_n$ to itself and some $n \geq 2$.  Conversely, 
all such algebras are weakly regular.
\end{theorem}
\noindent Zhang also examined the basic properties of such algebras. For $n > 2$, these are all non-noetherian domains of exponential growth which contain a free associative algebra, and so do not have a classical quotient division ring of fractions.  They have no normal elements.  

While for $n > 2$ these algebras are non-noetherian, Piontkovski proved that these algebras (and in fact all algebras defined by a single degree $2$ relation in a free algebra) are graded coherent---that is, every finitely generated graded left or right ideal is finitely presented \cite{Pi08}.  Coherence is a weak substitute for the noetherian property, which implies at least that the subcategory of finitely presented modules is abelian.  Just as it is conjectured that all regular algebras are noetherian, Bondal has conjectured that all weakly AS-regular algebras are graded coherent \cite[Conjecture 1.3]{Pi08}.  For $A$ as in Theorem~\ref{thm:zhang}, Piontkovski also studied the associated category $\rproj A$, which is defined here as the category of finitely presented graded $A$-modules modulo finite-dimensional modules.  These categories can be interpreted  as noncommutative projective lines.  Sisodia and Smith studied some further properties of these categories, including their Grothendieck groups \cite{SiSm15}.  Another point of view on the graded coherent property for these algebras (and more generally, preprojective algebras) was given by Minamoto \cite{Min18}.

As far as weakly regular algebras of dimension greater than 2, Smith gave an interesting such example of dimension $3$ with $7$ generators and $7$ relations, which is constructed from the geometry of the octonions \cite{Smi11}.  We mentioned earlier in Section~\ref{sec:extalgebra} the $m$-Koszul algebras defined by Berger.  Berger gives a family of examples defined by relations given by anti-symmetrizers \cite[Theorem 1.1]{BerMa06}, which include weakly AS-regular algebras which are $m$-Koszul, for unbounded choices of $m$.  Further examples of $m$-Koszul weakly regular algebras of dimension $3$ are the Yang-Mills algebras studied by Connes and Dubois-Violette \cite{CoDV02}.

Other examples have been studied in the more general context of twisted Calabi-Yau algebras, which we review in the next section.

\subsection{Twisted Calabi-Yau algebras}

Suppose that $A, B$ are any $\kk$-algebras (not necessarily graded).  We write $B^{op}$ for the opposite algebra of $B$.  Recall that an $(A, B)$-bimodule $M$ is the same as an $A \otimes_{\kk} B^{op}$-module, where $(a \otimes b) \cdot m = amb$ for $m \in M$, $a, b \in A$.  The algebra $A^e = A \otimes_{\kk} A$ is called the \emph{enveloping algebra} of $A$, and a left $A^e$-module is the same as an $(A, A)$-bimodule.  Because $(A^e)^{op} \cong A^{op} \otimes_k A \cong A \otimes_k A^{op}$ as algebras, we can also identify $(A, A)$-bimodules with right $A^e$-modules.  

\begin{definition}
Let $A$ be a $\kk$-algebra.  Suppose that 
\begin{enumerate}
\item $A$ is \emph{homologically smooth}:  that is, the $A^e$-module $A$ has a length $d$ projective resolution over $A^e$ by finitely generated projective left $A^e$-modules;  and
\item $\Ext^i_{A^e}(A, A^e) \cong \begin{cases} 0 & i \neq d \\ {} A^{\mu} & i = d \end{cases}$ as right $A^e$-modules, for some automorphism $\mu: A \to A$.
\end{enumerate}
then $A$ is called \emph{twisted Calabi-Yau}.
\end{definition}
Here, the module $A^{\mu}$ has its right action by $A$ twisted by the automorphism $\mu$, i.e. $b$ acts on the right of $a \in A^{\mu}$ by $a * b = a \mu(b)$.  The automorphism $\mu$ is called the \emph{Nakayama automorphism} of $A$.  It is determined only up to composition with an inner automorphism of $A$.  If one may take $\mu = 1$ in the definition, then one says that $A$ is \emph{Calabi-Yau}.  

Ginzburg defined the notion above and studied some properties of Calabi-Yau algebras in \cite{Gi}.   He was inspired by Kontsevich's notion of a Calabi-Yau category, which itself was related to investigations of mirror symmetry involving Calabi-Yau varieties in algebraic geometry.  Thus  Calabi and Yau ended up with their names on these algebraic objects via a rather roundabout route.  While there was reason to be particularly interested in Calabi-Yau algebras from this original motivation, from the algebraic perspective the generalization to twisted Calabi-Yau algebras is very natural, as we will see.

We now abbreviate Calabi-Yau by CY.  Note that the two parts of the definition of twisted CY algebra have clear formal similarities to the first and third axioms of the definition of AS regular algebra, where the special $(A, A)$-bimodule $A$ acts as a kind of generalized trivial module.  Indeed, algebras such as Sklyanin algebras were among Ginzburg's initial examples of CY algebras.  In \cite[Lemma 1.2]{RRZ14}, Reyes, the author, and Zhang proved that a connected $\mb{N}$-graded $\kk$-algebra is  weakly AS regular if and only if $A$ is twisted CY.  Technically, this result is about a version of the twisted CY property stated in the graded category, but this doesn't make any difference \cite[Theorem 4.2]{ReRo22}. One direction of \cite[Lemma 1.2]{RRZ14} is quite straightforward:  if $A$ is connected graded twisted CY, by tensoring the minimal $A^e$-projective resolution of $A$ with the trivial module $\kk$ on each side, one gets good properties for the minimal projective $A$-resolution of $\kk$, which implies that $A$ is weakly AS regular.  This direction was also shown in the special case of Calabi-Yau algebras in \cite[Proposition 4.3]{BeTa07}.  The converse uses the theory of dualizing complexes and relies on an earlier paper of Yekutieli and Zhang \cite{YeZh06}.  
%This theorem helps to explain why there is no standard term for  ``weakly AS regular algebra", as these days one can simply refer to these as connected graded twisted CY algebras.

From this point of view, we see that twisted CY algebras provide a vast generalization of  AS regular algebras to arbitrary, possibly ungraded algebras.  The main difference is the lack of a condition on growth in the twisted CY concept.  This simply reflects the fact that there was less attention initially from those interested in CY algebras on ring-theoretic properties such as the noetherian condition, which typically fails without a restriction on growth.

The notion of twisted CY unifies other possible generalizations of the (weakly) AS regular condition.
For example, suppose that $Q$ is a finite quiver with $n$ vertices, and let $\kk Q$ be its path algebra, graded by path length.  Let $A = \kk Q/I$ for some graded ideal $I$ generated by homogeneous elements of degree at least $2$.  Martinez-Villa and Solberg \cite{M-VSo11} defined a notion of generalized AS regular for such an algebra $A$, which essentially 
treats the $n$ simple left and right modules associated to the vertices of $Q$ as a set of trivial modules.  Namely, the AS Gorenstein condition is replaced by a requirement that $\Ext^j_A(-, A)$ acts as $0$ on all simple graded left modules for $0 \leq j < d = \gldim A$, while
$\Ext^d(-, A)$ gives a bijection from the simple graded left modules to the simple graded right modules (as a set).  It turns out that generalized Artin Schelter regular in this sense is the same as twisted CY for those $\mb{N}$-graded algebras of the form 
$\kk Q/I$ \cite[Theorem 3.19]{MiMo11}.  More generally, the author and Reyes characterized the finitely generated $\mb{N}$-graded twisted CY algebras $A$ with $\dim_{\kk} A_i < \infty$ over any base field $\kk$, in terms of another variation on the definition of AS regular \cite[Theorem 1.5]{ReRo22}. 

Twisted CY algebras are their own subject, with a large and growing literature we cannot hope to address further in this survey.  However, in the next section we describe some of the ways in which the language and ideas of CY algebras have played back into the theory of regular algebras in recent years.

\section{From twisted CY algebras back to AS regular algebras}
\label{sec:CYtoAS}

\subsection{Superpotentials}
Let $F = \kk \langle x_1, \dots x_m \rangle$ be a free algebra, graded with $\deg(x_i) = 1$ for all $i$, and suppose that $\tau$ is some graded automorphism of $F$.  Fixing a degree $n$, there is a unique $\kk$-linear map $\rho: F_n \to F_n$ which acts on the basis of monomials via $x_{i_1} \dots x_{i_n} \mapsto \tau(x_{i_n}) x_{i_1} \dots x_{i_{n-1}}$.  A \emph{$\tau$-twisted superpotential} is a homogeneous element $w \in F_n$ such that $\rho(w) = w$.  For each $1 \leq i \leq m$ there is an associated linear operator $\partial_i$ which acts on monomials in $F$ by $\partial_i(v) = u$ if $v = x_i u$ for some monomial $u$ and $\partial_i(v) = 0$ otherwise. Given some $k \leq n-2$, one may define an associated \emph{derivation quotient algebra} $D(w, k) = F/I$, where $I$ is generated by the set of  elements 
\[
\{ \partial_{i_1} \circ \partial_{i_2} \circ \dots \circ \partial_{i_k}(w) | 1 \leq i_1, \dots, i_k \leq m \}.
\]
This setup also generalizes easily to the quiver setting, by replacing $F$ with a path algebra $\kk Q$ and the variables with arrows.  Two overarching questions about such algebras are (i) whether every twisted CY algebra is a derivation quotient algebra, in graded or local settings; and (ii) for which $w, k$ the algebra $D(w, k)$ is twisted CY of dimension $k+2$, in which case the corresponding potential $w$ is called \emph{good} for the given $k$.

Let us first address progress on question (i).  Note that Artin and Schelter's original work already showed that regular algebras of dimension $3$ come from twisted superpotentials (just not in that language).  More exactly, as we explained in Section~\ref{sec:reg3}, any regular algebra of dimension $3$ either is equal to $D(w, 1)$ for some twisted potential $w \in \kk \langle x_1, x_2, x_3 \rangle_3$; or else equal to $D(w, 1)$ for some twisted potential $w \in \kk \langle x_1, x_2 \rangle_4$.  The restrictions on the number of generators and the degree of the potential are due to the polynomial growth condition.  

Bocklandt showed that Calabi-Yau algebras of dimension $3$ of the form $\kk Q/I$ are derivation quotient algebras $D(w, 1)$ \cite{Bock08}.  Then Bocklandt, Schedler and Wemyss \cite{BSW10} showed that any $m$-Koszul graded twisted CY algebra of the form $\kk Q/I$ with global dimension $d \geq 2$ is a derivation quotient algebra of the form $D(w, d-2)$.  Finally, Van den Bergh proved a very general result on the existence of superpotentials for complete CY algebras, which includes the graded case \cite{Vdb15}.

Question (ii) is much more sensitive.  It was already hard for regular algebras of dimension $3$, where it amounts to determining which values of the parameters in each family of the classification yield a regular algebra.  The deep geometric methods of \cite{ATV1} were introduced to answer this.     So the question about which potentials are good in an even more general setting is expected to be very difficult. 

In any case, the language and point of view of superpotentials and the associated concepts from CY algebras, in particular the Nakayama automorphism, have come to be increasingly important in the theory of regular algebras.  An early example is the paper \cite{D-V07}, in which Dubois-Violette reinterprets some of the theory of weakly regular algebras of dimension $2$ and $3$ using the language of potentials and the techniques of $m$-Koszul algebras. 

\subsection{Superpotentials of regular algebras}

The details of the classification of regular algebras of dimension $3$ have been reexamined in recent years, with more of a focus on the twisted superpotentials.  First, the parametrized superpotentials listed explicitly in \cite{AS87} do not produce all regular algebras of dimension $3$ up to isomorphism; they are ``generic forms", in the sense that they are restricted to those whose associated twist $\tau$ is diagonalizable.  Because diagonalizable matrices form an open subset of the general linear group, the superpotentials of other regular algebras lie in the closures of the orbits of the given superpotentials under the general linear group.  Second, although which twisted superpotentials are good can be checked, in theory, from the geometric results in \cite{ATV1}, a more direct understanding in terms of the properties of the superpotential is desirable.  

Mori and Smith considered the quadratic Calabi-Yau algebras of dimension $3$ in \cite{MoSm17}.  Given a superpotential $w$ of degree $3$ with associated derivation quotient algebra $A(w) = D(w, 1)$, the image of $w$ in the symmetric algebra has some zero locus in $\mb{P}^2$. A complete classification of which $w$ are good is given, primarily in terms of the geometry of that zero locus (which is different from the point scheme of $A(w)$).  
%One interesting consequence of this classification is that when $A(w)$ is not regular, it is not a domain and in fact there are two degree $1$ elements whose product is $0$.  

Mori and Ueyama extended similar methods to cubic Calabi-Yau algebras of dimension $3$ \cite{MoUe19}.  Extensions to twisted superpotentials, therefore recovering the classification of general AS regular algebras from the point of view of superpotentials, followed in several papers of Itaba and Matsuno  \cite{ItMa21} \cite{ItMa22}.  In related work, Iyudu and Shkarin \cite{IySh20} gave a classification of the isomorphism classes of the derivation quotient algebras $D(w, 1)$ for the $w$'s occurring in this setting (good or not), using more purely computational methods.   

In a different direction, Chirvasitu, Kanda, and Smith gave a new take on constructing normal extensions of regular algebras using twisted superpotentials, which works equally well in the setting of weakly regular algebras \cite{CKS19}.

\subsection{The Nakayama automorphism and the CY case}

The Nakayama automorphism $\mu$ of a regular algebra has received new attention as an important invariant.   In general, given a regular algebra defined by generators and relations, it may be hard to calculate the Nakayama automorphism. There are number of papers which show how to calculate $\mu_A$ for certain classes of algebras, see for instance \cite{LWW14}, \cite{ZVZ17}, \cite{LMZ17b}.

In the $m$-Koszul case, the Nakayama automorphism is closely related to the superpotential.  Suppose that $A$ is an AS regular derivation quotient algebra $A = D(w, k)$, for a $\mu$-twisted superpotential $w$. Here $\mu$ was defined as an automorphism of the free algebra, but one may see that by construction it preserves the relations of $A$, and so descends to an automorphism of $A$.  Importantly, this is essentially the same as the Nakayama automorphism $\mu$ of $A$ \cite[Theorem 4.4]{MoSm16}.  (Technically, this is true only up to sign, and possibly up to replacing $\mu$ by $\mu^{-1}$ depending on one's conventions.   There are sign conventions which motivated the ``super" part of the term superpotential, which we have ignored for simplicity).   

Notably, some of the regular algebras that have long attracted the most interest, such as Sklyanin algebras, are CY (where $\mu = 1$), and so there is reason to focus on this case in particular. An important example is the work of Pym, who studied the flat deformations of the polynomial ring $\kk[x_1, \dots, x_4]$---which automatically stay CY---and examined the resulting AS-regular algebras, which fall into $6$ families \cite{Pym15}.   Any such algebra gives a Poisson structure to the polynomial ring as a semi-classical limit, and the proof relies heavily on Poisson geometry.  Vancliff had studied earlier a more limited family of flat deformations coming from twists \cite{Van99}.  In related work, Lecoutre and Sierra generalized one of Pym's 6 families to give examples of CY regular algebras in all dimensions \cite{LeSi19}.

One would like to relate an arbitrary twisted CY algebra to a CY one in some way.  Reyes, the author, and Zhang studied the behavior of the Nakayama automorphism under Zhang twist \cite{RRZ14}. Suppose that $\kk$ is algebraically closed of characteristic $0$.  In this case, for a noetherian AS regular $A$ with Nakayama automorphism $\mu$  there is a always a graded automorphism $\sigma$ such that the associated Zhang twist $A^{\sigma}$ has Nakayama automorphism of finite order \cite[Theorem 5.3]{RRZ17}, \cite[Corollary 0.7]{RRZ14}.   Moreover, if $\mu$ has finite order, then the skew group algebra $A \rtimes \langle \mu \rangle$ is CY  of the same global dimension as $A$ (but is no longer connected graded) \cite[Theorem 5.3]{RRZ17}, \cite[Corollary 0.6]{RRZ14}. 

Alternatively, if one is willing to bump the global dimension up by one, the Ore extension $A[x; \mu]$ is always CY \cite[Theorem 5.3]{RRZ17},\cite[Corollary 0.6]{RRZ14}.  This was proved independently in the Koszul case in \cite{HVZ13b}, and Goodman and Kr\"ahmer proved a more general result for arbitrary (ungraded) algebras $A$ \cite{GoKr14}.    %Similarly, the fixed ring $A^{\langle \mu \rangle}$ has  trivial Nakayama automorphism but in general is only AS-Gorenstein, \cite[Corollary 0.6]{RRZ14}.
These results provide a potential framework for reducing the amount of work needed to prove that regular algebras have good properties.  For example, if all CY regular algebras of dimension $d+1$ are domains, then all regular algebras of dimension $d$ are domains.

We close the paper with a reminder that the origin of the subject---the paper \cite{AS87} by Artin and Schelter---proceeded by considering twisted superpotentials and classifying geometrically the generic forms of these that could lead to a regular algebra of dimension $3$ (though not in that language).  A direct analog of this for $4$-dimensional  regular algebras is still lacking.  There is not a good understanding of the parameter space of twisted superpotentials and how the known classes of $4$-dimensional regular algebras fit into it and relate to each other.  Some initial work on generalizing \cite{AS87} to dimension $4$ was done by Caines in his PhD thesis \cite{Cai03}, but it was just a start.  This idea seems ripe for further development.

\end{document}